\documentstyle{amsppt}
\nologo

\refstyle{A}

\hoffset .1 true in \voffset .2 true in

\hsize=6.3 true in \vsize=8.7 true in

\topmatter

\title Gap Sheaves and Vogel Cycles \endtitle

\author David B. Massey \endauthor

\address{David B. Massey, Dept. of Mathematics, Northeastern University,  Boston, MA, 02115, USA}
\endaddress

\email{DMASSEY\@neu.edu}\endemail



\endtopmatter

\document

\noindent\S0. {\bf Introduction}.  

\vskip .2in

Throughout our work on the L\^e cycles of an affine hypersurface singularity (see [{\bf M2-5}]), our primary algebraic tool consisted of a
method for taking the Jacobian ideal of a complex analytic function and decomposing it into pure-dimensional
``pieces''.  These pieces were obtained by considering the relative polar varieties of L\^e and Teissier (see, for example, [{\bf L-T}],
[{\bf T1}], [{\bf T2}]) as gap sheaves in the sense of [{\bf S-T}]. A gap sheaf is a formal device which gives a scheme-theoretic meaning to the analytic closure of the difference of an
initial scheme and an analytic set. 
 We would like to extend our results on L\^e cycles to functions on an
arbitrary complex analytic space, and so we need to generalize this algebraic approach.

We begin with an ordered set of generators for an ideal, and produce a collection of pure-dimensional analytic
cycles, the {\bf Vogel cycles} (see [{\bf G1}], [{\bf G2}], and [{\bf V}]), which seem to contain a great deal of ``geometric'' data related
to the original ideal. 
	The Vogel cycles are defined using gap sheaves, together with the associated analytic cycles which they define, the gap
cycles. 

 If the underlying analytic space is not Cohen-Macaulay, the main technical problem is that there are, at least, three different
reasonable definitions of the gap sheaves and cycles; we select as ``the'' definition the one that works most nicely in
inductive proofs. We show, however, that if one re-chooses the functions defining the ideal in a suitably ``generic'' way,
then all competing definitions for the gap cycles and Vogel cycles agree.

In Section 3, we prove some extremely general L\^e-Iomdine-Vogel formulas; these formulas generalize the L\^e-Iomdine formulas that we used
so profitably in [{\bf M2-5}].

\vskip .4in

A brief summary of our primary definitions and results follows.

\vskip .2in

 Throughout this introduction and the entire paper, if $A$ and $B$ are
analytic cycles which intersect properly inside an ambient manifold $M$, we let $A\cdot B$ denote the associated intersection product cycle
({\bf not} cycle class); see [{\bf Fu}].

\vskip .2in

Let  $W$  be analytic subset of an analytic space $X$  and let  $\alpha$  be a coherent sheaf of ideals in  ${\Cal
O}_{{}_X}$. Let $V$ denote the scheme $V(\alpha)$. Then, the {\it gap sheaf} $V\lnot W$ is the analytic closure of
$V-W$; that is, $V\lnot W$ is the scheme obtained from $V$ by removing any components or embedded subvarieties
contained in $W$.

\vskip .1in

Let $X$ be a $d$-dimensional irreducible (though not necessarily reduced) analytic space and
let $\bold f:= (f_0, \dots, f_k)\in \left(\Cal O_{{}_X}\right)^{k+1}$. {The \it $i$-th  gap variety of $\bold
f$}, $\Pi^i_{\bold f}$, is defined as 
$$
\Pi^i_{\bold f}:= V(f_{i+k+1-d}, \dots, f_k)\ \lnot\ V(\bold f),
$$
if $d-(k+1)< i<d$. Similarly, the {\it $i$-th  modified gap variety of $\bold f$},
$\widetilde\Pi^i_{\bold f}$, is defined as 
$$
\widetilde\Pi^i_{\bold f}:= V(f_{i+k+1-d}, \dots, f_k)\ \lnot\ V(f_{i+k-d}),
$$ if $d-(k+1)< i<d$.  The {\it $i$-th  inductive gap variety of $\bold f$}, $\widehat\Pi^i_{\bold f}$, is
defined by downward induction 

$$
\widehat\Pi^d_\bold f=\cases X, &\text{if } \bold f\not\equiv 0\\
\emptyset, &\text{if }\bold f\equiv 0 \endcases
$$
and
$$
\widehat\Pi^i_{\bold f}:= \left(\widehat\Pi^{i+1}_{\bold f}\cap V(f_{i+k+1-d})\right)\ \lnot\ V(f_{i+k-d}),
$$ if $d-(k+1)< i<d$.

\vskip .2in

If $X$ is irreducible and Cohen-Macaulay, and each $\Pi^i_{\bold f}$ is purely $i$-dimensional,then all three types of gap
varieties are equal. If $X$ is an arbitrary irreducible space, then, locally, we may replace each member of the
tuple $\bold f$ by a ``generic'' linear combination of the elements of $\bold f$ to obtain a new tuple, a {\it generic
linear reorganization of $\bold f$}, for which the gap sheaves, modified gap sheaves, and inductive gap sheaves are all
equal.

\vskip .1in

If $X$ is irreducible of dimension $d$ and each $\Pi^i_{\bold f}$ is purely $i$-dimensional, then, on $X-V(\bold f)$,
$\big[\widehat\Pi^i_{\bold f}\big] = V(f_k)\cdot V(f_{k-1})\cdot\dots\cdot V(f_{i+k+1-d})$; hence, on $X$,
$\big[\widehat\Pi^i_{\bold f}\big]$ is the closure in $X$ of this cycle on $X-V(\bold f)$. 

\vskip .1in

If $X$ is a union of irreducible components, $X= \bigcup_j X_j$, then we do not define gap {\bf sheaves}, but only gap
{\bf cycles}. Writing $\big[V\big]$ for the cycle defined by a scheme $V$, we define the
{\it $i$-th  gap cycle of $\bold f$} by $\left[\Pi^i_{\bold f}\right] :=
\sum_j\big[\Pi^i_{\bold f_{|_{X_j}}}\big]$,  the {\it $i$-th  modified gap cycle of $\bold f$} by
$\big[\widetilde\Pi^i_{\bold f}\big] :=
\sum_j\big[\widetilde\Pi^i_{\bold f_{|_{X_j}}}\big]$, and the {\it $i$-th  inductive gap cycle of $\bold f$} by
$\big[\widehat\Pi^i_{\bold f}\big] :=
\sum_j\big[\widehat\Pi^i_{\bold f_{|_{X_j}}}\big]$. 

More generally, if we have an analytic cycle $M:=\sum_l m_l [V_l]$ in $X$, where all of the $m_l$ have the same sign, then we define the
various gap cycles relative to
$M$ by taking the sum of the appropriate gap cycles restricted to each of the $V_l$, weighted by the $m_l$. The requirement that all of the
$m_l$ have the same sign prevents the cancellation of contributions from the various $V_l$.

\vskip .1in

The modified gap varieties and cycles are merely an intermediate tool. The inductive gap varieties are what we actually use
to define (below) our primary objects of study: the Vogel cycles. However, the hypotheses that must be satisfied before we
can define the Vogel cycles include, crucially, the hypothesis that each gap set $\big|\Pi^i_{\bold f}\big|$ has dimension
$i$. Thus, while one can safely forget the definition of the modified gap varieties, both the gap varieties and
inductive gap varieties are important for our future results.

\vskip .1in  

If $X$ is irreducible of dimension $d$ and each $\Pi^i_{\bold f}$ is purely $i$-dimensional, then the {\it $i$-th Vogel cycle
of $\bold f$}, $\Delta^i_\bold f$ is given by $$\Delta^i_\bold f=\Big(\big[\widehat\Pi^{i+1}_{\bold f}\big]\cdot
\big[V(f_{i+k+1-d})\big]\Big)- \big[\widehat\Pi^{i}_{\bold f}\big].$$ If $X$ is a
union of irreducible components, then the $i$-th Vogel cycle is obtained by summing the $i$-th Vogel cycles of all of the
irreducible components (as in the definition of the gap cycles). Similarly, one obtains Vogel cycles with respect to a given cycle $M$ by
taking the weighted sum of the Vogel cycles of $\bold f$ restricted to each subvariety appearing in $M$.

If each $\Delta^i_\bold f$ is purely $i$-dimensional (which one can obtain locally by replacing $\bold f$ by a generic linear
reorganization), then each Vogel cycle,
$\Delta^i_\bold f$, is non-negative  and is contained in $V(\bold f)$. Moreover, $V(\bold
f)=\bigcup_i \big|\Delta^i_\bold f\big|$. Thus, we think of the Vogel cycles as decomposing $V(\bold f)$ on the level of
cycles.

\vskip .2in

We proved the important {\it Segre-Vogel Relation}:  Let $X$ be an irreducible, $d$-dimensional, analytic subset of
an analytic manifold $\Cal U$, let
$\bold f = (f_0, \dots, f_k)\in \left(\Cal O_{{}_X}\right)^{k+1}$, let $\pi:\operatorname{Bl}_{\bold f}X\rightarrow X$
denote the blow-up of
$X$ along $\bold f$, and let $E_{\bold f}$ denote the corresponding exceptional divisor.

\vskip .1in

If $E_{\bold f}$ properly intersects $\Cal U\times\Bbb P^m\times\{\bold 0\}$ in $\Cal U\times\Bbb P^k$ for all $m$, then Vogel
cycles of
$\bold f$ are defined and, in a neighborhood of $V(\bold f)$, for all $i$,
$$\widehat\Pi^{i+1}_{\bold f} = \pi_*(\operatorname{Bl}_{\bold f}X\cdot (\Cal U\times \Bbb P^{i+k+1-d}\times\{\bold 0\}))$$ and
$$\Delta^i_{\bold f} = \pi_*(E_{\bold f}\cdot (\Cal U\times \Bbb P^{i+k+1-d}\times\{\bold 0\})),$$ where the intersection takes
place in
$\Cal U\times \Bbb P^k$ and $\pi_*$ denotes the proper push-forward.

\vskip .1in

Moreover, for all $\bold p\in X$, there exists an open neighborhood $\Cal U$ of $\bold p$ in $\Cal U$ such that, for a generic
linear reorganization,
$\tilde{\bold f}$, of
$\bold f$, 
$E_{\tilde{\bold f}}$ properly intersects
$\Cal U\times\Bbb P^m\times\{\bold 0\}$ inside
$\Cal U\times \Bbb P^k$ for all $m$. In the algebraic category, we may produce such generic linear reorganizations globally,
i.e., such that $E_{\tilde{\bold f}}$ properly intersects
$\Cal U\times\Bbb P^m\times\{\bold 0\}$ inside
$\Cal U\times \Bbb P^k$ for all $m$.

\vskip .1in

What we have just stated is the Segre-Vogel Relation for an irreducible space $X$, as it appears in Theorem 2.18. We give a more general
version with respect to a pure-dimensional cycle in Corollary 2.20.

\vskip .2in

Finally, we derive the L\^e-Iomdine-Vogel (LIV) formulas: Let $X$ be an irreducible analytic space of dimension
$n+1$, let $\bold f:= (f_0, \dots, f_n)\in
\big(\Cal O_X\big)^{n+1}$, let $g\in\Cal O_X$, and let $\bold p\in V(\bold f, g)$. Let $a$ be a non-zero complex number, let
$j\geqslant 1$ be an integer, and let $\bold h := (f_1, \dots, f_n, f_0+ag^j)$.

Suppose that the Vogel cycles of $\bold f$ are defined at $\bold p$, and that $V(g)$ properly intersects each of the Vogel
cycles, $\Delta^i_\bold f$, at
$\bold p$  for all $i\geqslant 1$.

\vskip .1in

If $j$ is sufficiently large, then there is an equality of sets given by
$V(\bold h) = V(\bold f, g)$,
${\operatorname{dim}_\bold p}V(\bold h) = \big({\operatorname{dim}_\bold p}V(\bold f)\big) -1$ provided that
${\operatorname{dim}_\bold p}V(\bold f)\geqslant 1$, the Vogel cycles of $\bold h$ exist at $\bold p$, and

\vskip .2in

$$\Delta^0_\bold h = \Delta^0_\bold f+j\big(\Delta^1_\bold f\cdot V(g)\big)$$ and, for $1\leqslant i\leqslant n-1$,
$$
\Delta^i_\bold h = j\big(\Delta^{i+1}_\bold f\cdot V(g)\big).
$$

\vskip .3in

In particular, if $j\geqslant 1+
\big(\Delta^0_\bold f\big)_\bold p$, then these conclusions hold. Once again, there is a more general version of this result with respect to
the cycle $M$.

\vskip .4in

\noindent\S1. {\bf Gap Sheaves}.  

\vskip .2in

 Let  $W$  be analytic subset of an analytic space $X$  and let  $\alpha$  be a coherent sheaf of ideals in  ${\Cal
O}_{{}_X}$.  At each point $\bold x$ of $V(\alpha)$, we wish to consider scheme-theoretically those components of 
$V(\alpha )$  which are not contained in $|W|$. 

\vskip .3in

\noindent{\bf Definition 1.1}. Let $A$ denote $\Cal O_{{}_{X, \bold x}}$; we  write $\alpha_\bold x$ for the stalk of
$\alpha$ in $A$. Let  $S$  be the multiplicatively closed set  $A - \bigcup {\slanted p}$   where the union is over all $
\slanted p\in\operatorname{Ass}(A/\alpha_\bold x) $ with   $|V({\slanted p})| \nsubseteq\ |W|$. Then, we define
$\alpha_\bold x\lnot W$  to equal  $S^{-1}\alpha_\bold x\cap A$.  Thus,  $\alpha_\bold x\lnot W$  is the ideal in $A$
consisting of the  intersection of those (possibly embedded) primary ideals, $\slanted q$,  associated to   
$\alpha_\bold x$  such that  $|V(q)|
\nsubseteq\ |W|$. 

Now, we have defined $\alpha_\bold x\lnot W$ in each stalk.  By [{\bf S-T}], if we perform this operation simultaneously
at all points of $V(\alpha)$, then we obtain a coherent sheaf of ideals called a {\it gap sheaf} ;  we write this sheaf as
$\alpha\lnot W$.  If  $V = V(\alpha)$, we let $V\lnot W$ denote  the scheme  $V(\alpha\lnot W)$.  

It is important to note that the scheme  $V\lnot W$  does not depend on the structure of  $W$  as a scheme, but only as an
analytic set.  The scheme $V\lnot W$ is sometimes referred to as the {\it analytic closure of $V-W$} [{\bf Fi}, p.41]; this
is certainly the correct, intuitive way to think of $V\lnot W$.

\vskip .05in

We find it convenient to extend this gap sheaf notation to the case of analytic sets (reduced schemes) and analytic cycles.  

Hence, if $Z$ and $W$ are analytic sets, then we let $Z\lnot W$ denote the union of the components of $Z$ which are not
contained in $W$; if $C=\sum m_i[V_i]$ is an analytic cycle in a complex manifold $M$ and $W$ is an analytic subset of $M$,
then we define $C\lnot W$ by 
$$ C\lnot W = \sum_{V_i\not\subseteq W} m_i[V_i] .
$$ If $\alpha$ is a coherent sheaf of ideals in $\Cal O_M$, $C$ is a cycle in $M$, and $W$ is an analytic subset of
$M$, then clearly
$[V(\alpha)\lnot W] = [V(\alpha)]\lnot W$ and $|C\lnot W|= |C|\lnot W$.

\vskip .3in

\noindent{\it Remark 1.2}. Later, the reader may wonder why we do not define something analogous to a gap sheaf, but where
we {\bf keep} those components which are contained in a given analytic set, $W$, instead of throwing them away. There are
several reasons for this.

First of all, on the level of schemes, we can not make this approach work; the primary ideals in a primary decomposition (of
a given ideal) which define varieties contained in $W$ would not be independent of the decomposition. We could just take the
isolated primary ideals which define varieties contained in $W$, but this disposes of too much algebraic structure.

Our best attempt at taking a scheme $V$ and keeping, scheme-theoretically, those (possibly embedded) components which are
contained in $W$ would be to consider $V\lnot(V\lnot W)$.

However, even this device would not aid us much later; as we shall see -- beginning with Definition 2.14 -- we need to deal
more with the intersection product on analytic cycles, and not so much with primary decompositions.

\vskip .3in

The following lemma is very useful for calculating $V\lnot W$.

\vskip .2in

\noindent{\bf Lemma 1.3}. {\it Let $(X, \Cal O_{{}_X})$ be an analytic space, let $\alpha, \beta$, and $\gamma$ be coherent
sheaves of ideals in  $\Cal O_{{}_X}$, let $f,g \in \Cal O_{{}_X}$, and let $W$, $Y$, and $Z$ be  analytic subsets of
$X$ such that $Z\subseteq W$.  Then,

\vskip .1in

\noindent i)\hskip .32in  $\alpha\lnot W = (\alpha\lnot Z)\lnot W$, and thus, as schemes, $V(\alpha)\lnot W =
(V(\alpha)\lnot Z)\lnot W$;

\vskip .1in

\noindent ii)\hskip .25in  $(\alpha + \beta)\lnot W \ = \ (\alpha\lnot Z + \beta)\lnot W$, and thus, as schemes,
 $$\big(V(\alpha) \cap V(\beta)\big)\lnot W\ =\ \big(V(\alpha\lnot Z) \cap V(\beta)\big)\lnot W ;$$

\vskip .1in

\noindent iii)\hskip .25in  if \ $V(\alpha + \gamma) \subseteq W$, then $\big((\alpha \cap \beta) + \gamma\big)\lnot W\ =\
(\beta +
\gamma)\lnot W$, and thus, as schemes,
  $$\Big(\big(V(\alpha) \cup V(\beta)\big) \cap V(\gamma)\Big)\lnot W\ = \ \big(V(\beta) \cap V(\gamma)\big)\lnot W ;$$

\vskip .1in

\noindent iv)\hskip .25in  if \ $V(\alpha + <g>) \subseteq W$, then $(\alpha + <fg>)\lnot W \ = \ (\alpha + <f>)\lnot W$, and
thus, as schemes, 
$$\big(V(\alpha) \cap V(fg)\big)\lnot W \ = \  \big(V(\alpha) \cap V(f)\big)\lnot W .$$

\vskip .1in

\noindent v)\hskip .27in  $\alpha\lnot (W\cup Y) =  (\alpha\lnot W)\lnot Y$, and
thus, as schemes, 
$$V(\alpha)\lnot (W\cup Y) =  (V(\alpha)\lnot W)\lnot Y.$$

\vskip .1in

\noindent The analog of ii) for sets and cycles is also trivial to verify; that is,
$$\big(|V(\alpha)| \cap |V(\beta)|\big)\lnot W\ =\ \Big(\big(|V(\alpha)|\lnot Z\big) \cap |V(\beta)|\Big)\lnot W,$$ and, if
all intersections are proper,
$$\big([V(\alpha)] \cdot [V(\beta)]\big)\lnot W\ =\ \Big(\big([V(\alpha)]\lnot Z\big) \cdot [V(\beta)]\Big)\lnot W.$$
 }

\vskip .3in

\noindent{\it Proof}. Statements i), ii), iii), and iv) are merely exercises in localization (see [{\bf M3}]). Statement v) is immediate.\qed

\vskip .4in

\noindent{\it Remark 1.4}. While it is a trivial observation, it is frequently important and useful to note that, for any
coherent sheaf of ideals,
$\alpha$,  in
$\Cal O_{{}_X}$ and for any $f\in \Cal O_{{}_X}$, $V(\alpha)\lnot V(f)$ and $V(f)$ intersect properly and $V(f)$ contains no
embedded subvarieties of $V(\alpha)\lnot V(f)$; thus, the intersection cycle  $[V(\alpha)\lnot V(f)]\cdot [V(f)]$ in $X$ is
well-defined (without having to mention an ambient manifold) and is equal to $$[V(<\alpha\ \lnot\ V(f)>+<f>)].$$

If $V(\alpha)$ and $V(f)$ intersect properly, then $[V(\alpha)]=[V(\alpha)\lnot V(f)]$ and, hence, 
$$ [V(\alpha)]\cdot[V(f)] = [V(\alpha)\lnot V(f)]\cdot[V(f)] = [V(<\alpha\ \lnot\ V(f)>+<f>)].
$$

\vskip .4in  

\noindent{\bf Lemma 1.5}. {\it Let $X$ be purely $d$-dimensional and Cohen-Macaulay. Let $f_1, \dots, f_k\in\Cal O_X$ and let
$W$ be an analytic subset of $X$.  If $V(f_1, \dots, f_k)\ \lnot\ W$ is purely $(d-k)$-dimensional, then it contains no
embedded subvarieties.

}

\vskip .3in

\noindent{\it Proof}. By definition,  $V(f_1, \dots, f_k)\ \lnot\ W$ can not have any embedded subvarieties contained in
$W$. At points, $\bold p$, outside of
$W$, $f_1, \dots, f_k$ determines a regular sequence in the Cohen-Macaulay ring $\Cal O_{X, \bold p}$ ; hence, there are no
embedded subvarieties outside of
$W$.\qed

\vskip .4in

\noindent{\it Example 1.6}. For the remainder of this section, we wish to describe the blow-up of a space along an ideal;
the description via gap sheaves is very nice. We shall use this description in the next section.

\vskip .1in

Let $(X, \Cal O_{{}_X})$ be an analytic space, and let $\bold f:=(f_0, \dots, f_k)$ be an ordered $(k+1)$-tuple of elements
of
$\Cal O_{{}_X}$. Then, the {\it blow-up of $X$ along $\bold f$} consists of an analytic subspace
$\operatorname{Bl}_{\bold f}X\subseteq X\times\Bbb P^k$, together with the projection morphism $\pi
:\operatorname{Bl}_{\bold f}X\rightarrow X$, which is the restriction of the standard projection from $X\times\Bbb P^k$ to
$X$.  If we use $[w_0: \dots : w_k]$ for homogeneous coordinates on $\Bbb P^k$, then the blow-up is given as a scheme by
$$
\operatorname{Bl}_{\bold f}X := V\big(\{w_if_j-w_jf_i\}_{{}_{\ 0\leqslant i, j\leqslant k}}\big) \ \lnot\  \left(V(f_0,
\dots, f_k)\times\Bbb P^k\right) .
$$

\vskip .2in

In order to describe the exceptional divisor as a cycle, we need to work on affine coordinate patches in $\Bbb P^k$.  We
shall  describe both the blow-up and the exceptional divisor on each affine patch $\{w_j\neq 0\}$.

On the patch $\{w_j\neq 0\}$, we use coordinates $\widetilde w_i := w_i/w_j$ for all $i\neq j$.
 Then, 
$$
\{w_j\neq 0\}\cap\operatorname{Bl}_{\bold f}X \ = \ V(\{f_i-\widetilde w_if_j\}_{{}_{i\neq j}})\ 
\lnot \left(V(f_j)\times\Bbb P^k\right) ,\tag{1.7}
$$  and the exceptional divisor, $E$, is the cycle defined on each affine patch in the following manner
$$
\{w_j\neq 0\}\cap E \ := \ \big[V\big(\{f_i-\widetilde w_if_j\}_{{}_{i\neq j}}\big)
\lnot \left(V(f_j)\times\Bbb P^k\right)\big]\cdot \left[V(f_j)\times\Bbb P^k\right] .\tag{1.8} 
$$ 

\vskip .1in

We have made these definitions with respect to a chosen $(k+1)$-tuple $\bold f$.  In fact, the analytic isomorphism-type of
the morphism $\pi :\operatorname{Bl}_{\bold f}X\rightarrow X$ only depends on the ideal, $I$, generated by the components
$f_0,
\dots, f_k$; this isomorphism-type is referred to as the {\it blow-up of $X$ along $I$}.  Of course, the isomorphism-type of
the exceptional divisor also depends only on the ideal $I$, and this isomorphism-type is simply called the {\it the
exceptional divisor of the blow-up of $X$ along $I$}.

\vskip .4in

\noindent  \S2. {\bf Gap Cycles and Vogel Cycles}

\vskip .2in 

Let $X$ be a $d$-dimensional analytic space and let $\bold f := (f_0, \dots, f_k)$ be an ordered $(k+1)$-tuple of elements of
$\Cal O_{{}_X}$.   We will define a sequence of cycles, the {\it Vogel cycles} ([{\bf Vo}], [{\bf Gas1}], [{\bf Gas2}]) of
$\bold f$; these cycles provide  effectively calculable data about the coherent sheaf of ideals $<f_0, \dots, f_k>$. Before
we can define the Vogel cycles, we must first define the {\it gap varieties} and {\it gap cycles} of $\bold f$.

It is useful to define gap and Vogel objects with respect to a given cycle  (see Section 4 for a discussion of {\bf why} the introduction of
a cycle is so useful). Hence, throughout this paper, we let
$M$ denote the cycle $\sum_l m_l[V_l]$ in $X$; we assume that this is a minimal presentation of $M$ -- that is, we assume that the
$V_l$ are distinct, irreducible, analytic subsets of $X$ and that none of the $m_l$ equal zero. In addition, to avoid cancellation of
contributions from various $V_l$, we assume that all of the $m_l$ have the same sign, i.e., that $\pm M>0$.

\vskip .2in

If $X$ is a union of irreducible components $\{X_i\}$, we will define the gap and Vogel cycles in $X$ as sums of the gap and
Vogel cycles from each $X_i$; similarly, we will define gap and Vogel cycles with respect to $M$ simply by taking weighted sums of
the gap and Vogel cycles of the irreducible components. The case of an irreducible space $X$ can be recovered from the cycle case by
simply taking $M=[X]$.   Thus, we find that we need to first define the gap varieties, gap cycles, and Vogel cycles in the case where $X$ is
irreducible.

However, even if we assume that the underlying space is irreducible, there is a further complication in the general setting:
$\Cal O_X$ may not be Cohen-Macaulay. This causes numerous problems, for we must worry about embedded subvarieties. To deal
with this problem, we introduce three avatars of gap varieties and examine the relations between them. 

We will define the
(ordinary) gap varieties, $\big\{\Pi^i_\bold f\big\}_i$, the modified gap varieties, $\big\{\widetilde\Pi^i_\bold f\big\}_i$,
and the inductive gap varieties, $\big\{\widehat\Pi^i_\bold f\big\}_i$. We shall use the inductive gap varieties to define the
Vogel cycles, but need to make assumptions about the (ordinary) gap varieties in order for the definition to make sense; the
modified gap varieties are merely a convenient tool for proving results about $\Pi^i_\bold f$ and 
$\widehat\Pi^i_\bold f$.

\vskip .2in

\noindent{\bf Definition 2.1}. Assume that $X$ is irreducible (though, not necessarily reduced). 
For all
$i$, we define the {\it gap varieties}, the {\it modified gap varieties}, and the {\it inductive gap varieties} of $\bold f$,
which we denote by
$\Pi^i_\bold f$,
$\widetilde\Pi^i_\bold f$, and $\widehat\Pi^i_\bold f$, respectively.

\vskip .1in

First, if $i<d-(k+1)$ or $i>d$, we set $\Pi^i_{\bold f} = \widetilde\Pi^i_{\bold f} = \widehat\Pi^i_{\bold f} = \emptyset$.

\vskip .2in

We define 
$\Pi^d_\bold f := X\lnot V(\bold f)$, and $\widetilde\Pi^d_\bold f = \widehat\Pi^d_\bold f:=X\lnot V(f_k)$.
\vskip .3in

For $d-(k+1)\leqslant i<d$,  the {\it $i$-th  gap variety of $\bold f$}, $\Pi^i_{\bold f}$, is defined as 
$$
\Pi^i_{\bold f}:= V(f_{i+k+1-d}, \dots, f_k)\ \lnot\ V(\bold f).
$$ Note that $\Pi^{d-(k+1)}_{\bold f}=\emptyset$.

\vskip .3in

For $d-(k+1)< i<d$, the {\it $i$-th  modified gap variety of $\bold f$}, $\widetilde\Pi^i_{\bold f}$, is defined as 
$$
\widetilde\Pi^i_{\bold f}:= V(f_{i+k+1-d}, \dots, f_k)\ \lnot\ V(f_{i+k-d});
$$ we define $\widetilde\Pi^{d-(k+1)}_{\bold f}:= \emptyset$.

\vskip .3in

For $d-(k+1)< i<d$, the {\it $i$-th  inductive gap variety of $\bold f$}, $\widehat\Pi^i_{\bold f}$, is defined by downward
induction (recall that  $\widehat\Pi^d_\bold f$ is defined above)
$$
\widehat\Pi^i_{\bold f}:= \left(\widehat\Pi^{i+1}_{\bold f}\cap V(f_{i+k+1-d})\right)\ \lnot\ V(f_{i+k-d});
$$ we define $\widehat\Pi^{d-(k+1)}_{\bold f}:= \emptyset$.

\vskip .2in

Naturally, we define the {\it $i$-th gap cycle}, {\it modified $i$-th gap cycle}, and {\it inductive $i$-th gap cycle} of
$\bold f$ to be the cycles defined by these schemes, i.e., 
$\big[\Pi^i_{\bold f}\big]$, $\big[\widetilde\Pi^i_{\bold f}\big]$, and $\big[\widehat\Pi^i_{\bold f}\big]$ , respectively.

\vskip .2in

If $X$ is a union of irreducible components $\{X_j\}$ and $\bold f\in \left(\Cal O_{{}_X}\right)^{k+1}$, then we define the
{\it
$i$-th  gap cycle of $\bold f$} by $\left[\Pi^i_{\bold f}\right] := \sum_j\big[\Pi^i_{\bold f_{|_{X_j}}}\big]$,  the {\it
$i$-th  modified gap cycle of $\bold f$} by $\big[\widetilde\Pi^i_{\bold f}\big] :=
\sum_j\big[\widetilde\Pi^i_{\bold f_{|_{X_j}}}\big]$, and the {\it $i$-th  inductive gap cycle of $\bold f$} by
$\big[\widehat\Pi^i_{\bold f}\big] :=
\sum_j\big[\widehat\Pi^i_{\bold f_{|_{X_j}}}\big]$. 

We define the {\it
$i$-th  gap set of $\bold f$}, the {\it $i$-th  modified gap set of $\bold f$}, and the {\it $i$-th  inductive gap set of
$\bold f$} to be $\Big|\big[\Pi^i_{\bold f}\big]\Big|$, 
$\Big|\big[\widetilde\Pi^i_{\bold f}\big]\Big|$, and $\Big|\big[\widehat\Pi^i_{\bold f}\big]\Big|$, respectively. We will write
simply $\big|\Pi^i_{\bold f}\big|$, 
$\big|\widetilde\Pi^i_{\bold f}\big|$, and $\big|\widehat\Pi^i_{\bold f}\big|$, respectively.

\vskip .1in

Finally, we need to define gap cycles and sets with respect to the cycle $M$. We define the
{\it
$i$-th  gap cycle of $\bold f$ with respect to $M$} by $\Pi^i_{\bold f}(M) := \sum_l m_l \big[\Pi^i_{\bold
f_{|_{V_l}}}\big]$,  the {\it
$i$-th  modified gap cycle of $\bold f$ with respect to $M$} by $\widetilde\Pi^i_{\bold f}(M) :=
\sum_l m_l\big[\widetilde\Pi^i_{\bold f_{|_{V_l}}}\big]$, and the {\it $i$-th  inductive gap cycle of $\bold f$ with respect to
$M$} by
$\widehat\Pi^i_{\bold f}(M) :=
\sum_l m_l \big[\widehat\Pi^i_{\bold f_{|_{V_l}}}\big]$.

Of course, we define the associated gap sets with respect to $M$ to be
the sets underlying the various gap cycles.

\vskip .2in

\noindent Note that we have not defined gap {\bf varieties} for $\bold f$ unless $X$ is irreducible.

\vskip .3in

The following proposition gives a number of basic results and interrelationships between the various gap varieties.

\vskip .2in

\noindent{\bf Proposition 2.2}. Let $X$ be irreducible. Then, 

\vskip .1in

\noindent i)\hskip .1in there is an inclusion of sets
$\big|\widehat\Pi^i_{\bold f}\big|\subseteq
\big |\widetilde\Pi^i_{\bold f}\big |\subseteq
\big |\Pi^i_{\bold f}\big |$, \ $\big|\widehat\Pi^i_{\bold f}\big|$ is purely $i$-dimensional, and all components\linebreak \hbox{}\hskip
.2inof 
$\big |\widetilde\Pi^i_{\bold f}\big |$ and $\big |\Pi^i_{\bold f}\big |$ have
dimension at least
$i$;

\vskip .1in

\noindent ii)\hskip .07in  there is an equality of schemes $\Pi^{i}_{\bold f} = \big(\Pi^{i+1}_{\bold f}\cap
V(f_{i+k+1-d})\big)\lnot V(\bold f)$; 

\vskip .1in

\noindent iii) if $i\leqslant d$, then $\big|\Pi^{i-1}_{\bold f}\big|\subseteq \big|\Pi^i_{\bold f}\big|$ and
$\big|\widehat\Pi^{i-1}_{\bold f}\big|\subseteq \big|\widehat\Pi^i_{\bold f}\big|$;

\vskip .1in

\noindent iv)\hskip .07in the sets $V(f_{i+k+1-d})$ and $\big|\widehat\Pi^{i+1}_{\bold f}\big|$ intersect properly, and
there is an equality of cycles
$$\big[\widehat\Pi^{i}_{\bold f}\big] =
\Big(\big[\widehat\Pi^{i+1}_{\bold f}\big]\cdot
\big[V(f_{i+k+1-d})\big]\Big)\ \lnot\ V(f_{i+k-d}) ;$$

\vskip .1in

\noindent v)\hskip .1in if there is an equality of sets $\big|\widetilde\Pi^i_{\bold f}\big| = \big|\Pi^i_{\bold f}\big|$,
then the schemes $\widetilde\Pi^i_{\bold f}$ and $\Pi^i_{\bold f}$ are equal up to embedded subvariety, and so there is an
equality of cycles 
$\big[\widetilde\Pi^i_{\bold f}\big] =
\big[\Pi^i_{\bold f}\big]$;

\vskip .1in

\noindent vi)\hskip .07in if there is an equality of sets $\big|\widehat\Pi^j_{\bold f}\big| = \big|\Pi^j_{\bold f}\big|$
for all
$j\geqslant i+1$, then
$\big[\widehat\Pi^i_{\bold f}\big] \leqslant
\big[\Pi^i_{\bold f}\big]$;

\vskip .1in

\noindent vii)\hskip .04in if there is an equality of schemes $\Pi^{i+1}_{\bold f} = \widehat\Pi^{i+1}_{\bold f}$, then
there is an equality of schemes $\widetilde\Pi^{i}_{\bold f} = \widehat\Pi^{i}_{\bold f}$.

\vskip .3in

\noindent{\it Proof}. i) is obvious from the definitions.  ii) follows immediately from Lemma 1.3.ii (using $V(\bold f)$ for
both
$Z$ and
$W$).  v) is immediate.

\vskip .1in

\noindent Proof of iii): ii) implies $\big|\Pi^{i-1}_{\bold f}\big|\subseteq \big|\Pi^i_{\bold f}\big|$. That
$\big|\widehat\Pi^{i-1}_{\bold f}\big|\subseteq \big|\widehat\Pi^i_{\bold f}\big|$ follows from the inductive definition. 

\vskip .1in

\noindent Proof of iv):  By definition, $\widehat\Pi^{i+1}_{\bold f}$ has no components or embedded subvarieties contained in
$V(f_{i+k+1-d})$.  Thus, $\big[\widehat\Pi^{i+1}_{\bold f}\cap V(f_{i+k+1-d})\big] = \big[\widehat\Pi^{i+1}_{\bold
f}\big]\cdot
\big[V(f_{i+k+1-d})\big]$.  The desired conclusion follows.

\vskip .1in

\noindent Proof of vi): By downward induction on $i$.  Note first that $\big[\widehat\Pi^d_{\bold f}\big] \leqslant
\big[\Pi^d_{\bold f}\big]$, since they are, in fact, equal.  Suppose now that $i<d$ and that there is an equality of sets
$\big|\widehat\Pi^j_{\bold f}\big| = \big|\Pi^j_{\bold f}\big|$ for all $j\geqslant i+1$.  From induction, we know that
$\big[\widehat\Pi^{i+1}_{\bold f}\big] \leqslant
\big[\Pi^{i+1}_{\bold f}\big]$.   Thus,
$$
\big[\Pi^i_{\bold f}\big] = \big[\Pi^{i+1}_{\bold f}\cap V(f_{i+k+1-d})\big]\ \lnot\ V(\bold f)\ \geqslant\
\big[\Pi^{i+1}_{\bold f}\cap V(f_{i+k+1-d})\big]\ \lnot\ V(f_{i+k-d}).
$$

\vskip .1in

Since iv) tells us that $\big|\Pi^{i+1}_{\bold f}\big|$ intersects $V(f_{i+k+1-d})$ properly, we may apply [{\bf Fu}, 8.2.a]
to conclude that 
$\big[\Pi^{i+1}_{\bold f}\cap V(f_{i+k+1-d})\big]\geqslant \big[\Pi^{i+1}_{\bold f}\big]\cdot
\big[V(f_{i+k+1-d})\big]$ (the presence of embedded varieties in $\Pi^{i+1}_{\bold f}$ can cause a strict inequality).
Therefore, 
$$
\big[\Pi^i_{\bold f}\big]\geqslant \Big(\big[\Pi^{i+1}_{\bold f}\big]\cdot
\big[V(f_{i+k+1-d})\big]\Big)\ \lnot\ V(f_{i+k-d}).
$$ Now, applying our inductive hypothesis and iv), we conclude that $\big[\widehat\Pi^i_{\bold f}\big] \leqslant
\big[\Pi^i_{\bold f}\big]$.

\vskip .1in

\noindent Proof of vii): We have
$$
\widehat\Pi^i_{\bold f} = \big(\widehat\Pi^{i+1}_{\bold f}\cap V(f_{i+k+1-d})\big)\ \lnot\ V(f_{i+k-d}) =
\big(\Pi^{i+1}_{\bold f}\cap V(f_{i+k+1-d})\big)\ \lnot\ V(f_{i+k-d}).
$$ By 1.3.i, this equals $\Big(\big(\Pi^{i+1}_{\bold f}\cap V(f_{i+k+1-d})\big)\ \lnot\ V(\bold f)\Big)\ \lnot\
V(f_{i+k-d})$.  By ii) of this proposition, this last expression equals $\Pi^i_{\bold f}\ \lnot\ V(f_{i+k-d}) =
\big(V(f_{i+k+1-d}, \dots, f_k)\
\lnot\ V(\bold f)\big)\ \lnot\ V(f_{i+k-d})$.  Applying 1.3.i again, we find that $\widehat\Pi^i_{\bold f} = V(f_{i+k+1-d},
\dots, f_k)\ \lnot\ V(f_{i+k-d})=\widetilde\Pi^i_{\bold f}$.\qed

\vskip .5in

We wish to define the Vogel cycles now.  However, before we can do this, we need to decide which of the different gap cycles
to use to define the Vogel cycles.  As a preliminary step, we first define the sets which will underlie the Vogel cycles.

\vskip .2in

\noindent{\bf Definition 2.3}. Assume that $X$ is irreducible. If $i\neq d$, then we define the {\it $i$-th Vogel set of
$\bold f$}, $D^i_{\bold f}$, to be the union of the irreducible components of $\big|\Pi^{i+1}_{\bold f}\cap
V(f_{i+k+1-d})\big|$ which are contained in $V(\bold f)$; by 2.2.ii, this is equivalent to
$$ D^i_{\bold f}= \overline{\big|\Pi^{i+1}_{\bold f}\cap V(f_{i+k+1-d})\big|\ -\ \big|\Pi^{i}_{\bold f}\big|} .
$$

\vskip .2in

\noindent We set $D^d_{\bold f}=\cases \emptyset, &\text{if }\bold f\not\equiv 0\\X, &\text{if }\bold f\equiv 0.\endcases$

\vskip .2in

\noindent Note that, if $i<d-(k+1)$ or $i>d$, then $D^i_{\bold f}=\emptyset$.

\vskip .2in

\noindent If $X$ is a union of irreducible components $\{X_j\}$, we define $D^i_{\bold f}:= \bigcup_j D^i_{{\bold
f}_{|_{X_j}}}$.

\vskip .2in

\noindent We define the {\it $i$-th Vogel set of $\bold f$ with respect to $M$} to be  $D^i_{\bold f}(M):= \bigcup_{j} D^i_{{\bold
f}_{|_{V_l}}} .$

\vskip .4in

\noindent{\bf Proposition 2.4}.  {\it Every component of $D^i_{\bold f}(M)$ has dimension at least $i$ and $|M|\cap V(\bold f) =
\bigcup_i D^i_{\bold f}(M)$. If $\Pi^i_\bold f(M)$ is $i$-dimensional and $C$ is an $i$-dimensional irreducible component of
$|M|\cap |V(\bold f)|$, then
$C\subseteq D^i_\bold f(M)$. 

If $X$ is irreducible of dimension $d$, then, for all $i\leqslant d-1$,
$$\big|\Pi^{i+1}_{\bold f}\cap V(\bold f)\big|=
\bigcup_{k\leqslant i}D^k_{\bold f}.$$

}

\vskip .3in

\noindent{\it Proof}.   We may work on each irreducible set, $V_l$, separately; therefore, we assume that we are in
the case where $X$ is irreducible and $M=[X]$. 

\vskip .1in

That every component of $D^i_{\bold f}$ has dimension at least $i$ follows immediately from the fact that each component of
$\Pi^{i+1}_{\bold f}$ has dimension at least $i+1$ (by 2.2.i).

\vskip .1in

By definition $X = \big|\Pi^d_{\bold f}\big| \cup D^d_{\bold f}$. Hence, $V(\bold f) = \big|\Pi^d_{\bold f}\cap V(\bold
f)\big|
\cup D^d_{\bold f}$ and so, the equation $\dsize V(\bold f) = \bigcup_i D^i_{\bold f}$  follows once we show the final claim
of the proposition.

\vskip .1in

Suppose that  $i\leqslant d-1$. Then, 
$$\big|\Pi^{i+1}_{\bold f}\cap V(\bold f)\big|= \big|\Pi^{i+1}_{\bold f}\cap V(f_{i+k+1-d})\cap V(\bold f)\big|=
\big|\big(\Pi^i_{\bold f}\cup D^i_{\bold f}\big)\cap V(\bold f)\big| =  \big|\Pi^{i}_{\bold f}\cap V(\bold f)\big|\cup
D^i_{\bold f}.$$ As $\Pi^{i}_{\bold f}$ is eventually empty, the desired conclusion follows.

\vskip .1in

Finally, suppose that $C$ is an $i$-dimensional irreducible component of $|V(\bold f)|$ and $\Pi^i_\bold f$ is
$i$-dimen\-sional. Then, $C$ is contained in a component $C^\prime$ of $|V(f_{i+k+2-d}, \dots, f_k)|$; such a
$C^\prime$ necessarily has dimension at least $i+1$. Thus,
$C^\prime$ cannot be contained in $V(\bold f)$. It follows that $C^\prime$ is contained in $\Pi^{i+1}_\bold f$. Therefore, 
$$ C\ \subseteq\ C^\prime\cap V(f_{i+k+1-d})\ \subseteq\ \big|\Pi^{i+1}_\bold f\cap V(f_{i+k+1-d})\big|\ =\
\big|\Pi^i_\bold f\big|\cup D^i_\bold f .
$$ If $\Pi^i_\bold f$ is $i$-dimensional, then -- since $C\subseteq V(\bold f)$ and is $i$-dimensional -- it follows that
$C\not\subseteq
\big|\Pi^i_\bold f\big|$, and so $C\subseteq D^i_\bold f$.
\qed

\vskip .3in

Below, we prove the {\it Dimensionality Lemma} in which we state as hypotheses/conclusions that ``$\big|\Pi^i_{\bold
f}(M)\big|$ is purely $i$-dimensional'' and ``$D^i_{\bold f}(M)$ is purely $i$-dimensional''. Since sets cannot be
negative-dimensional, for $i<0$, we mean that the respective set is empty.  

\vskip .3in

\noindent{\bf Lemma 2.5} ({\bf Dimensionality Lemma}). {\it The following are equivalent:

\vskip .1in

\noindent i)\hskip .28in for all $i$, $\big|\Pi^i_{\bold f}(M)\big|$ is purely $i$-dimensional;

\vskip .1in

\noindent ii)\hskip .24in  for all $i$, $\big|\Pi^i_{\bold f}(M)\big| = \big|\widetilde\Pi^i_{\bold f}(M)\big|$;

\vskip .1in

\noindent iii)\hskip .2in  for all $i$, $\big|\Pi^i_{\bold f}(M)\big| = \big|\widehat\Pi^i_{\bold f}(M)\big|$.

\vskip .2in

In addition, these equivalent conditions imply 

\vskip .1in

\noindent iv)\hskip .2in for all $i$, $D^i_{\bold f}(M)$ is purely
$i$-dimensional;

\vskip .1in

\noindent and, for all $\bold p\in |M|\cap V(\bold f)$, there exists a neighborhood of $\bold p$ in which iv) implies i), ii), and
iii). }

\vskip .3in

\noindent{\it Proof}.  Again we may consider each component appearing $M$ separately; hence, we may assume that $X$ is
irreducible and $M=[X]$.

\vskip .2in

 As all the statements are set-theoretic, to cut down on notation, we shall omit the vertical lines around the various gap
sheaves. 

\vskip .1in

We will show that i) and iii) are each equivalent to ii), that i) implies iv), and that, near points of $V(\bold f)$, iv)
implies i).

\vskip .2in

\noindent i)$\Rightarrow$ ii):\hskip .2in  Assume i).  From the definition of $\Pi^i_{\bold f}$, what we need to show is: if
$C$ is a component of $V(f_{i+k+1-d}, \dots, f_k)$, then $C$ is contained in $V(\bold f)$ if and only if $C$ is contained in
$V(f_{i+k-d})$.  As $V(\bold f)\subseteq V(f_{i+k-d})$, one implication is trivial, and so what we must show is that if 
$C$ is a component of $V(f_{i+k+1-d}, \dots, f_k)$ and $C\subseteq V(f_{i+k-d})$, then $C\subseteq V(\bold f)$.  

Suppose not. As $C$ is a component of $V(f_{i+k+1-d}, \dots, f_k)$, the dimension of $C$ is at least $i$. If
$C\not\subseteq V(\bold f)$, then -- by definition -- $C$ is a component of $\Pi^i_{\bold f}$.  But $C$ is also contained in
$V(f_{i+k-d})$, and so
$C$ is a component of $\Pi^i_{\bold f}\cap V(f_{i+k-d}) = \Pi^{i-1}_{\bold f}\cup D^{i-1}_{\bold f}$.  As $C$ is not
contained in
$V(\bold f)$, we conclude that $C$ is a component of $\Pi^{i-1}_{\bold f}$ of dimension at least $i$.  This contradicts i).

\vskip .2in

\noindent ii)$\Rightarrow$ i):\hskip .2in Assume ii).  From Definition 2.1, $\Pi^i_{\bold f}$ is purely
$i$-dimensional for $i\geqslant d$.  Suppose that $i_0$ is the largest integer $i$ (less than $d$) such that
$\Pi^i_{\bold f}$ is not purely $i$-dimensional.  Then, $\Pi^{i_0+1}_{\bold f}$ is purely $(i_0+1)$-dimensional and, by ii),
the set $\Pi^{i_0+1}_{\bold f}$ is equal to $V(f_{i_0+k+2-d}, \dots, f_k)\ \lnot \ V(f_{i_0+k+1-d})$.   Hence, the
intersection
$\Pi^{i_0+1}_{\bold f}\cap V(f_{i_0+k+1-d})$ is proper, and so $\Pi^{i_0+1}_{\bold f}\cap V(f_{i_0+k+1-d})$ is purely
$i_0$-dimensional.  As there is an equality of sets $\Pi^{i_0+1}_{\bold f}\cap V(f_{i_0+k+1-d}) = \Pi^{i_0}_{\bold f}\cup
D^{i_0}_{\bold f}$, this contradicts the fact that $\Pi^{i_0}_{\bold f}$ is not purely $i_0$-dimensional. 

\vskip .2in

\noindent iii)$\Rightarrow$ ii):\hskip .2in Assume iii). Then ii) follows immediately from the fact that 
$\widehat\Pi^i_{\bold f}\subseteq
\widetilde\Pi^i_{\bold f}\subseteq
\Pi^i_{\bold f}$ (see 2.2.i).

\vskip .2in

\noindent ii)$\Rightarrow$ iii):\hskip .2in Assume ii). The proof is by induction. iii) is certainly true by definition for
$i\geqslant d$. Now, suppose that $\Pi^i_{\bold f} =\widehat\Pi^i_{\bold f}$ for $i\geqslant m$, where $m\leqslant d$.  We
need to show that $\Pi^{m-1}_{\bold f} =
\widehat\Pi^{m-1}_{\bold f}$. We have
$$
\widehat\Pi^{m-1}_{\bold f} = \big(\widehat\Pi^{m}_{\bold f}\cap V(f_{m+k-d})\big)\ \lnot\ V(f_{m+k-1-d}) =
\big(\Pi^{m}_{\bold f}\cap V(f_{m+k-d})\big)\ \lnot\ V(f_{m+k-1-d}).$$
 By combining the definition of $\Pi^{m}_{\bold f}$ as $V(f_{m+k+1-d}, \dots, f_k)\
\lnot\ V(\bold f)$ with Lemma 1.3.ii, we conclude that 
$$\big(\Pi^{m}_{\bold f}\cap V(f_{m+k-d})\big)\ \lnot\ V(f_{m+k-1-d}) = V(f_{m+k-d}, \dots, f_k)\ \lnot\ V(f_{m+k-1-d})$$
and so,
$\widehat\Pi^{m-1}_{\bold f} = \widetilde\Pi^{m-1}_{\bold f}$. By ii), this implies that $\widehat\Pi^{m-1}_{\bold f} =
\Pi^{m-1}_{\bold f}$ and we are finished.

\vskip .2in

\noindent i)$\Rightarrow$ iv):\hskip .2in Assume i), and suppose that $i_0$ is such that $D^{i_0}_{\bold f}$ is {\bf not}
purely
$i_0$-dimensional.  Then, $\Pi^{i_0+1}_{\bold f}\cap V(f_{i_0+k+1-d})$ is not purely $i_0$-dimensional.  As 
$\Pi^{i_0+1}_{\bold f}$ is purely $(i_0+1)$-dimensional by assumption, it follows that $V(f_{i_0+k+1-d})$ contains a
component,
$C$, of $\Pi^{i_0+1}_{\bold f}$.  As $C$  is a component of $\Pi^{i_0+1}_{\bold f}$, $C$ is not contained in
$V(\bold f)$. 

Thus,  $C$ is a component of $\Pi^{i_0+1}_{\bold f}\cap V(f_{i_0+k+1-d}) = \Pi^{i_0}_{\bold f}\cup D^{i_0}_{\bold f}$ which
is not contained in $V(\bold f)$, and so $C$ is an $(i_0+1)$-dimensional component of 
$\Pi^{i_0}_{\bold f}$ -- this contradicts our assumption.

\vskip .2in

\noindent iv)$\Rightarrow$ i):\hskip .2in  Assume iv), and that we are interested in the germ of the situation at a point
$\bold p\in V(\bold f)$. Let
$i_0$ be the smallest
$i$ such that
$\Pi^i_{\bold f}$ is {\bf not} purely $i$-dimensional.  By Proposition 2.2.i,  $\Pi^{i_0}_{\bold f}$ must have dimension at
least
$i_0+1$. Thus, since $\bold p\in V(\bold f)$,
$\Pi^{i_0}_{\bold f}\cap V(f_{i_0+k-d})$ has dimension at least $i_0$.  But, as sets, $\Pi^{i_0}_{\bold f}\cap
V(f_{i_0+k-d}) =
\Pi^{i_0-1}_{\bold f}\cup D^{i_0-1}_{\bold f}$, and by assumption $D^{i_0-1}_{\bold f}$ is purely
$(i_0-1)$-dimensional.  Therefore, we conclude that $\Pi^{i_0-1}_{\bold f}$ has dimension at least $i_0$ -- a contradiction
of the choice of $i_0$.
\qed

\vskip .4in

\noindent{\it Remark 2.6}. Our phrasing of Lemma 2.5 is the most elegant, and is in the form that we will usually need.
However, it is occasionally helpful to note that our proof does not require that one knows i), ii), or iii) for {\bf all}
$i$. Specifically, what our proof actually shows is that:

\vskip .1in

\noindent $\bullet$ \hskip .2in if $\big|\Pi^{i-1}_{\bold f}(M)\big|$ is purely $(i-1)$-dimensional, then
$\big|\Pi^i_{\bold f}(M)\big| =
\big|\widetilde\Pi^i_{\bold f}(M)\big|$; 

\vskip .2in

\noindent $\bullet$ \hskip .2in if $\big|\Pi^i_{\bold f}(M)\big| = \big|\widetilde\Pi^i_{\bold f}(M)\big|$ for all $i>k$, then
$\big|\Pi^{k}_{\bold f}(M)\big|$ is purely $k$-dimensional;

\vskip .2in

\noindent $\bullet$ \hskip .2in if $\big|\Pi^i_{\bold f}(M)\big| = \big|\widehat\Pi^i_{\bold f}(M)\big|$, then
$\big|\Pi^i_{\bold f}(M)\big| = \big|\widetilde\Pi^i_{\bold f}(M)\big|$;

\vskip .2in

\noindent $\bullet$ \hskip .2in if $\big|\Pi^i_{\bold f}(M)\big| = \big|\widehat\Pi^i_{\bold f}(M)\big|$ for all
$i\geqslant m$, and
 $\big|\Pi^{m-1}_{\bold f}(M)\big| = \big|\widetilde\Pi^{m-1}_{\bold f}(M)\big|$, then $\big|\Pi^{m-1}_{\bold f}(M)\big|
=$\linebreak\hbox{}\hskip .3in
$\big|\widehat\Pi^{m-1}_{\bold f}(M)\big|$; in particular, if $\big|\Pi^i_{\bold f}(M)\big| =
\big|\widetilde\Pi^i_{\bold f}(M)\big|$ for all $i\geqslant m-1$, then $\big|\Pi^i_{\bold f}(M)\big| =
\big|\widehat\Pi^i_{\bold f}(M)\big|$\linebreak\hbox{}\hskip .3in for all $i\geqslant m-1$;

\vskip .2in

\noindent $\bullet$ \hskip .2in if $\big|\Pi^i_{\bold f}(M)\big|$ is purely $i$-dimensional and $\big|\Pi^{i+1}_{\bold f}(M)\big|$
is purely $(i+1)$-dimensional, then $D^i_\bold f(M)$ is \linebreak\hbox{}\hskip .3in purely $i$-dimensional; and

\vskip .2in

\noindent $\bullet$ \hskip .2in if $\bold p\in |M|\cap V(\bold f)$, $\big|\Pi^{i-1}_{\bold f}(M)\big|$ is purely
$(i-1)$-dimensional at $\bold p$, and $D^{i-1}_{\bold f}(M)$ is purely $(i-1)$-\linebreak\hbox{}\hskip .3in  dimensional
at $\bold p$, then $\big|\Pi^{i}_{\bold f}(M)\big|$ is purely $i$-dimensional at $\bold p$.

\vskip .4in

\noindent{\bf Definition 2.7}. If the equivalent conditions i), ii), and iii) of Lemma 2.5 hold, we say that {\it the gap
sets of $\bold f$ with respect to $M$ have the correct dimension}.

If the equivalent conditions i), ii), iii), and iv) of Lemma 2.5 hold at a point $\bold p\in |M|\cap V(\bold f)$, we say that {\it
the Vogel sets of $\bold f$ with respect to $M$ have the correct dimension at $\bold p$}. We say simply that {\it the Vogel sets of
$\bold f$ with respect to $M$ have the correct dimension} provided that they have correct dimension at all points of $|M|\cap
V(\bold f)$.

\vskip .3in

\noindent{\it Remark 2.8}. Note that, since  every component of $D^i_{\bold f}(M)$ has dimension at least $i$ (see 2.4), if the
Vogel sets all have correct dimension at $\bold p$, then all the Vogel sets have correct dimension at points {\bf near}
$\bold p$.

Note also that if the gap sets have correct dimension, then the Vogel sets have correct dimension. Moreover, since we are 
interested only in what happens near $V(\bold f)$, the natural assumption for us to make seems like it should be that the
Vogel cycles have correct dimension. However, our usual assumption will be that {\bf gap sets} have correct dimension; for
2.5 tells us that, in a neighborhood of $V(\bold f)$, these assumptions are equivalent, and requiring the gap sets to have
the correct dimension saves us from having to state over and over again that we take a small neighborhood of a point of
$V(\bold f)$.

It is important to remember that one
implication of the Vogel and gap sets having correct dimension is that
$D^i_{\bold f}(M)$,
$\Pi^{i}_{\bold f}(M)$, $\widetilde\Pi^{i}_{\bold f}(M)$, and $\widehat\Pi^{i}_{\bold f}(M)$ are all empty if
$i< 0$, and $\Pi^{0}_{\bold f}(M) = \widetilde\Pi^{0}_{\bold f}(M)=\widehat\Pi^{0}_{\bold f}(M) = \emptyset$ at points of
$|M|\cap V(\bold f)$. 

Finally, consider the special case where $\bold p$ is an isolated point of $|M|\cap V(\bold f)$. Then, 2.4 implies that, near $\bold p$,
$D^i_{\bold f}(M) =\emptyset$ if $i\geqslant 1$, and $D^0_{\bold f}(M)=\{\bold p\}$. Thus, 2.5 implies that the gap sets and the Vogel sets
have correct dimension at $\bold p$.

\vskip .4in

\noindent{\bf Proposition 2.9}. {\it If $X$ is irreducible and Cohen-Macaulay, and all of the gap sets of $\bold f$ have
correct dimension, then, for all
$i$, the schemes
$\Pi^i_{\bold f}$,
$\widetilde\Pi^i_{\bold f}$, and $\widehat\Pi^i_{\bold f}$ are equal. }

\vskip .3in

\noindent{\it Proof}.  By 2.2.v, if the gap sets have correct dimension, then the schemes $\Pi^i_{\bold f}$ and
$\widetilde\Pi^i_{\bold f}$ are equal up to embedded subvariety.  By Lemma 1.5, $\Pi^i_{\bold f}$ and
$\widetilde\Pi^i_{\bold f}$ have no embedded subvarieties; therefore, they are equal as schemes.

To prove that the scheme structure of $\widehat\Pi^i_{\bold f}$ agrees with the other two, we must, of course, use
induction.  Let
$d$ denote the dimension of $X$.  For $i\geqslant d$, we know that $\Pi^i_{\bold f}=\widetilde\Pi^i_{\bold
f}=\widehat\Pi^i_{\bold f}$. 

Suppose, inductively, that
$\Pi^{i+1}_{\bold f}=\widetilde\Pi^{i+1}_{\bold f}=\widehat\Pi^{i+1}_{\bold f}$. Then, 2.2.vii tells us that
$\widetilde\Pi^{i}_{\bold f} =
\widehat\Pi^{i}_{\bold f}$ and, by the first paragraph above, we know that this equals $\Pi^{i}_{\bold f}$. \qed

\vskip .3in

While we have been selecting $(k+1)$-tuples, $\bold f$, our primary object of interest is, in fact, the ideal $<\bold f>$
generated by the $f_0, \dots, f_k$. As far as the ideal $<\bold f>$ is concerned, the functions comprising $\bold f$ may not
be suitably generic. However, as we shall see, to obtain a well-behaved ordered collection of generators, one only needs to
replace
$(f_0,
\dots, f_k)$ by generic linear combinations of the $f_i$'s themselves. However, the term ``generic'' here is used in a
non-standard way; what we need is to replace $f_0$ by a generic linear combination, then -- fixing this new $f_0$ -- replace
$f_1$ by a generic linear combination, and so on. Since ``generic'' should always mean open and dense in {\bf some} topology,
we will define a new, convenient one.

\vskip .3in

\noindent{\bf Definition 2.10}. The {\it pseudo-Zariski topology} (pZ-topology) on a topological space $(X, \Cal T)$ is a new
topological space $(X, \Cal T_{pZ})$ given by $\Cal U\in\Cal T_{pZ}$ if and only if $\Cal U$ is empty or is an open, dense
subset in $(X, \Cal T)$. (One verifies easily that this, in fact, yields a topology on $X$.)

Given two topological spaces $X$ and $Y$, let $\pi_X$ and $\pi_Y$ denote the projections from $X\times Y$ onto $X$ and $Y$,
respectively. The {\it inductive pseudo-Zariski topology} (IPZ-topology) on $X\times Y$ is given by: $\Cal W\subseteq X\times
Y$ is open in the IPZ-topology if and only if $\pi_X(\Cal W)$ is open in the pZ-topology on $X$ and, for all $x\in \pi_X(\Cal
W)$,
$\pi_Y\big(\Cal W\cap \pi_X^{-1}(x)\big)$ is open in the pZ-topology on $Y$. (It is trivial to verify that this is a topology
on $X\times Y$, and that a non-empty open set in the IPZ-topology on $X\times Y$ is a dense set in the cross-product topology
on $X\times Y$.)

Finally, given a finite number of topological spaces $X_1$, $X_2$, $\dots$, $X_m$, the IPZ-topology on
$X_1\times X_2\times\dots\times X_m$ is given inductively by using the IPZ-topology on each product in the expression
$\Big(\big((X_1\times X_2)\times X_3\big)\times\dots\times X_{m-1}\Big)\times X_m$.

\vskip .1in

 A {\it generic linear reorganization} of a $(k+1)$-tuple $\bold f$ is a matrix product  $\bold fA$, where the matrix $A$ is
invertible and is an element of some given generic subset in the IPZ-topology on the $(k+1)$-fold product $\Bbb
C^{k+1}\times\dots\times\Bbb C^{k+1}$ (where we consider each column of $A$ to be contained in one copy of $\Bbb C^{k+1}$).

\vskip .3in

Note that, if $X_1=X_2=\dots = X_m=\Bbb C^N$ (or $\Bbb P^N$), then the IPZ-topology on the product is more fine than the
Zariski topology, but sets which are open in the IPZ-topology need {\bf not} be open in the classical topology on the
product.

\vskip .3in

\noindent{\bf Proposition/Definition 2.11}. {\it If $X$ is irreducible, then, for all
$\bold p\in X$, for a generic linear reorganization,
$\hat{\bold f}$, of
$\bold f$, the gap sets of $\hat{\bold f}$  all have correct dimension at $\bold p$ and, for all $i$, there is
an equality of schemes
$\Pi^i_{\hat{\bold f}}=\widetilde\Pi^i_{\hat{\bold f}}=\widehat\Pi^i_{\hat{\bold f}}$ in a neighborhood of $\bold p$.

Therefore, for $\bold p\in |M|$, for a generic linear
reorganization,
$\hat{\bold f}$, of
$\bold f$, the gap sets of $\hat{\bold f}$ with respect to $|M|$ all have correct dimension at $\bold p$ and, for all $i$, there is
an equality of cycles
$\Pi^i_{\hat{\bold f}}(M)=\widetilde\Pi^i_{\hat{\bold f}}(M)=\widehat\Pi^i_{\hat{\bold f}}(M)$ at $\bold
p$.}

\vskip .1in

If we are working in the algebraic category, then we may produce such generic linear reorganizations globally.

\vskip .1in

We refer to a reorganization $\hat{\bold f}$ such that the above equality of cycles holds as an {\it agreeable
reorganization of
$\bold f$ (with respect to $M$ at $\bold p$)} (for it makes the various cycle structures agree).

\vskip .3in

\noindent{\it Proof}. Assume that $X$ is irreducible. We fix a point $\bold p\in X$. Our sole reason for
stating the results ``at $\bold p$'' is that, at several places in the proof, we will need to know that certain analytic sets
have a finite number of analytic components. This is, of course, guaranteed near a given point or in the algebraic category.
Hence, throughout the proof, we will make no further reference to working in a neighborhood of $\bold p$, but will assume
that all of the analytic sets that arise have a finite number of components.

\vskip .1in

We first show:

\vskip .1in

\noindent $(\dagger)$ for a generic linear reorganization, $\hat{\bold f}$, of $\bold f$, for all $i$, $V(\hat f_{i+k-d})$
contains no component or embedded subvariety of $\Pi^i_{\hat{\bold f}}$.

\vskip .1in

We produce the $(k+1)$-tuple $\hat{\bold f}$ one element at a time, by downward induction. If $\bold f$ is identically zero
on
$X$, then
$(\dagger)$ is trivial. So, suppose that one of the $f_i$ does not vanish on $X$. Then, for a generic linear combination
$\hat f_k := a_0f_0 +\dots + a_kf_k$,
$\hat f_k$ does not vanish on $X$. Thus, $V(\hat f_{k})$ contains no component or embedded subvariety of
$\Pi^d_{\hat{\bold f}}$.

Now, suppose that we have made generic linear reorganizations of $\bold f$ to produce $\hat{\bold f}$, and that $V(\hat
f_{i+k-d})$  contains no component or embedded subvariety of $\Pi^i_{\hat{\bold f}}$ for all $i\geqslant m$. Then, for every
component or embedded subvariety, $W$, of $\Pi^{m-1}_{\hat{\bold f}}$, $W$ is contained in $V(\hat f_{m+k-d},
\dots, \hat f_k)$, but there exists some $\hat f_j$ with $j<m+k-d$ such that $W\not\subseteq V(\hat f_j)$. Thus, a generic
linear combination of the $\hat f$'s will not vanish on any component of embedded subvariety of
$\Pi^{m-1}_{\hat{\bold f}}$. This proves $(\dagger)$.

\vskip .1in

As $\Pi^{i-1}_{\hat{\bold f}} = \big(\Pi^i_{\hat{\bold f}}\cap V(\hat f_{i+k-d})\big)\lnot V(\hat{\bold f})$ by 2.2.ii,
$(\dagger)$ implies that the Vogel sets of $\hat{\bold f}$ all have correct dimension. We show that $\Pi^i_{\hat{\bold
f}}=\widetilde\Pi^i_{\hat{\bold f}}=\widehat\Pi^i_{\hat{\bold f}}$ by downward induction on $i$. 

When $i=d$, the statement is clear. Assume now that $\Pi^{i+1}_{\hat{\bold f}}=\widetilde\Pi^{i+1}_{\hat{\bold
f}}=\widehat\Pi^{i+1}_{\hat{\bold f}}$. Then,
$$
\widehat\Pi^i_{\hat{\bold f}} = \big(\widehat\Pi^{i+1}_{\hat{\bold f}}\cap V(\hat f_{i+k+1-d})\big)\lnot V(\hat f_{i+k-d}) =
\big(\Pi^{i+1}_{\hat{\bold f}}\cap V(\hat f_{i+k+1-d})\big)\lnot V(\hat f_{i+k-d}),
$$ which, by 1.3.i, equals $\Big(\big(\Pi^{i+1}_{\hat{\bold f}}\cap V(\hat f_{i+k+1-d})\big)\lnot V(\hat{\bold f})\Big)\lnot
V(\hat f_{i+k-d})$. 

Therefore, applying 2.2.ii, followed by $(\dagger)$, we conclude that $\widehat\Pi^i_{\hat{\bold f}} =
\Pi^i_{\hat{\bold f}}\
\lnot\  V(\hat f_{i+k-d}) = \Pi^i_{\hat{\bold f}}$.

\vskip .1in

As $\widehat\Pi^i_{\hat{\bold f}} = \Pi^i_{\hat{\bold f}}$ for all $i$, by applying  2.2.vii, we conclude that
$\widetilde\Pi^i_{\hat{\bold f}} =\widehat\Pi^i_{\hat{\bold f}} = 
\Pi^i_{\hat{\bold f}}$.\qed

\vskip .5in

We now wish to endow the Vogel sets a cycle structure. First, we need the following easy proposition.

\vskip .3in

\noindent{\bf Proposition 2.12}. {\it If $X$ is irreducible and the gap sets of $\bold f$ have correct
dimension, then there is an equality of cycles given by 
$$
\big[\widehat\Pi^{i}_{\bold f}\big] =\ \Big(\big[\widehat\Pi^{i+1}_{\bold f}\big]\cdot
\big[V(f_{i+k+1-d})\big]\Big)\ \lnot\ V(\bold f).
$$

}

\vskip .3in

\noindent{\it Proof}.  Note that, by 2.2.iv, the statement $\big[\widehat\Pi^{i}_{\bold f}\big] =\
\Big(\big[\widehat\Pi^{i+1}_{\bold f}\big]\cdot
\big[V(f_{i+k+1-d})\big]\Big)\ \lnot\ V(\bold f)$ is equivalent to the set-theoretic statement $
\big|\widehat\Pi^{i}_{\bold f}\big| =\ \Big(\big|\widehat\Pi^{i+1}_{\bold f}\big|\cap V(f_{i+k+1-d})\Big)\ \lnot\ V(\bold f)
$.  This set-theoretic statement follows easily from 2.5.iii; for it tells us that 
$$\Big(\big|\widehat\Pi^{i+1}_{\bold f}\big|\cap V(f_{i+k+1-d})\Big)\ \lnot\ V(\bold f) = \Big(\big|\Pi^{i+1}_{\bold
f}\big|\cap V(f_{i+k+1-d})\Big)\ \lnot\ V(\bold f)$$ and   2.2.ii tells us that this equals $\big|\Pi^i_{\bold f}\big|$. 
Applying 2.5.iii again yields the desired equality of cycles. \qed

\vskip .4in

\noindent{\it Remark 2.13}. It is  important to note that 2.12 gives a rather simple characterization of
$\big[\widehat\Pi^i_{\bold f}\big]$ near points where the Vogel sets have correct dimension: if $X$ is irreducible of
dimension
$d$, then, on $X-V(\bold f)$, $\big[\widehat\Pi^i_{\bold f}\big] = V(f_k)\cdot V(f_{k-1})\cdot\dots\cdot V(f_{i+k+1-d})$;
hence, on $X$,
$\big[\widehat\Pi^i_{\bold f}\big]$ is the closure in $X$ of this cycle on $X-V(\bold f)$. This cycle structure turns out to
be the correct one to use in order to endow the Vogel sets with a cycle structure. 

However, in order to guarantee that the
cycles we define actually have as their underlying sets the Vogel sets of $\bold f$, we only define the Vogel {\bf cycles}
when the Vogel sets have the correct dimensions and, even then, we must restrict ourselves to what happens in a neighborhood
of $V(\bold f)$.

\vskip .3in

\noindent{\bf Definition 2.14}. If $X$ is irreducible, and the Vogel sets of $\bold f$ all have correct
dimension, then we define the {\it $i$-th Vogel cycle of $\bold f$}, $\Delta^i_\bold f$, to be the sum of the components of 
$$
\Big(\big[\widehat\Pi^{i+1}_{\bold f}\big]\cdot
\big[V(f_{i+k+1-d})\big]\Big)- \big[\widehat\Pi^{i}_{\bold f}\big]
$$ which intersect $V(\bold f)$. In other words, if 
$$
\Big(\big[\widehat\Pi^{i+1}_{\bold f}\big]\cdot
\big[V(f_{i+k+1-d})\big]\Big)- \big[\widehat\Pi^{i}_{\bold f}\big]\ =\ \sum_j p_j\big[W_j\big],
$$ then $\dsize\Delta^i_\bold f\  = \ \sum_{W_j\cap V(\bold f)\neq\emptyset}p_j\big[W_j\big]$.

\vskip .1in

If all the Vogel sets of $\bold f$ with respect to $M$ have correct dimension, then we
say that {\it the Vogel cycles of $\bold f$ with respect to $M$ are defined} and their definition is  $\dsize\Delta^i_{\bold f}(M)
:=
\sum_l m_l\ 
\Delta^i_{{\bold f}_{|_{V_l}}}$.

\vskip .4in

Note that there is no difference between saying that the Vogel sets have correct dimension and that the Vogel cycles are
defined; we prefer to say that the Vogel cycles are defined, as the Vogel cycles are the objects in which we are most
interested. Notice, also, that Proposition 2.12 implies that, if the gap sets of $\bold f$ with respect to $M$ have correct
dimension, then
$\Delta^i_\bold f(M)$ is given simply by
\hbox{$\big(\widehat\Pi^{i+1}_{\bold f}(M)\cdot
V(f_{i+k+1-d})\big)- \widehat\Pi^{i}_{\bold f}(M)$}.

\vskip .4in

\noindent{\bf Proposition 2.15}. {\it If the Vogel cycles of $\bold f$ with respect to $M$ are defined, then each $\Delta^i_{{\bold
f}_{|_{V_l}}}$ is non-negative and purely $i$-dimensional.  Moreover, $\big|\Delta^i_{\bold f}(M)\big| = D^i_{\bold f}(M)\subseteq
|M|\cap V(\bold f)$.

}

\vskip .3in

\noindent{\it Proof}. As usual, we may assume that $X$ is irreducible and that $M=[X]$. That $\Delta^i_{\bold f}$ is
non-negative follows from the definition of the inductive gap varieties. That
$\Delta^i_{\bold f}$ is purely $i$-dimensional follows from the Dimensionality Lemma and the fact that 
$\widehat\Pi^{i+1}_{\bold f}$ and  $V(f_{i+k+1-d})$ intersect properly.  That $\big|\Delta^i_{\bold f}\big| = D^i_{\bold f}$
follows from the definition of $D^i_{\bold f}$ and the fact that the Dimensionality Lemma tells us that the gap varieties
and the inductive gap varieties are equal as sets in a neighborhood of $V(\bold f)$.\qed

\vskip .4in

\noindent{\it Remark 2.16}. If $X$ is irreducible, Proposition 2.12 and Proposition 2.15, together with the Dimensionality Lemma, tell us how
the Vogel cycles should be calculated; we will describe this now, omitting the square brackets for the cycles.

\vskip .1in

\noindent One begins with $\widehat\Pi^d_{\bold f} = X\lnot V(f_k)$; thus,  $\widehat\Pi^d_{\bold f}$  is either $0$ or $X$. Next,  one
calculates the intersection $\widehat\Pi^d_{\bold f}\cdot V(f_k)$. This intersection cycle has components contained in $V(\bold f)$ and
components which are not contained in
$V(\bold f)$. By 2.12, the sum of the components which are not contained in $V(\bold f)$ is precisely
$\widehat\Pi^{d-1}_{\bold f}$ and the sum of the components which are contained in $V(\bold f)$ is $\Delta^{d-1}_{\bold f}$.
Having calculated $\widehat\Pi^d_{\bold f}\cdot V(f_k) = \widehat\Pi^{d-1}_{\bold f} + \Delta^{d-1}_{\bold f}$, we use our
newly found
$\widehat\Pi^{d-1}_{\bold f}$ in the next step: the calculation of $\widehat\Pi^{d-1}_{\bold f}\cdot V(f_{k-1})$. One
proceeds downward inductively.

	The subtle point in the above description is that, if one is working in a neighborhood of a point of $V(\bold f)$, one may check {\bf while}
performing the calculation that the Vogel sets,
$\big|\Delta^i_{\bold f}\big|$, have correct dimension. For, by splitting the intersections into pieces which are, and
pieces which are not, contained in $V(\bold f)$, we are actually obtaining a cycle $\Delta^i_{\bold f}$ whose underlying set
is precisely
$D^i_{\bold f}$ (this follows from 2.2.ii).  Thus, one proceeds with the inductive calculation described above, and then
checks that the calculated $\Delta^i_{\bold f}$ have correct dimension, which then tells one that the calculation is
actually correct.

\vskip .1in

Consider the special case where $\bold p$ is an isolated point of $|M|\cap V(\bold f)$. As we saw in Remark 2.8, it is automatic that the
Vogel cycles are defined at $\bold p$, and only $\Delta^0_{\bold f}$ can be non-zero. In fact, if $X$ is irreducible of dimension $d$, then
Remark 2.13 implies that $\Delta^0_{\bold f}= V(f_k)\cdot\dots\cdot V(f_{k+1-d})$.

\vskip .4in

\noindent{\it Example 2.17}. We continue to suppress the square brackets around cycles.  Let $X= \Bbb C^5$ and let
$$\bold f = (f_0, f_1, f_2, f_3, f_4) =(-2ux^2, \ -2vx^2, \ -2wx^2, \ -3x^2 - 2x(u^2 + v^2 + w^2), \ 2y).$$ (The reason for
the strange, seemingly pointless, coefficients is that we will use this example later in a different context. See Example II.2.4.)
Then, $V(\bold f) = V(x,y)$ and $\widehat\Pi^5_{\bold f} = \Bbb C^5$.

$$\widehat\Pi^5_{\bold f}\cdot V(f_4) = \widehat\Pi^5_{\bold f}\cdot V(-2y) = V(y) .$$

As $V(y)$ is not contained in $V(\bold f)$, $\widehat\Pi^4_{\bold f} = V(y)$, and we continue.

\vskip .1in

$$\widehat\Pi^4_{\bold f} \cdot V(f_3) = V(y) \cdot V(-3x^2 - 2x(u^2 + v^2 + w^2)) = V(-3x - 2(u^2 + v^2 + w^2), y) + V(x,y)
=
\widehat\Pi^3_{\bold f} + \Delta^3_{\bold f} .$$

\vskip .1in

$$\widehat\Pi^3_{\bold f} \cdot V(f_2) =  V(-3x - 2(u^2 + v^2 + w^2), y) \cdot V(-2wx^2) =$$ 
$$V(-3x - 2(u^2 + v^2), w, y) + 2V(u^2 + v^2 + w^2, x, y) = \widehat\Pi^2_{\bold f} + \Delta^2_{\bold f} .$$

\vskip .1in

$$\widehat\Pi^2_{\bold f} \cdot V(f_1) = V(-3x - 2(u^2 + v^2), w, y) \cdot V(-2vx^2) = $$
$$V(-3x-2u^2, v, w, y) + 2V(u^2 + v^2, w, x, y) = \widehat\Pi^1_{\bold f} + \Delta^1_{\bold f} .$$

\vskip .1in

$$\widehat\Pi^1_{\bold f} \cdot V(f_0) = V(-3x-2u^2, v, w, y) \cdot V(-2ux^2) = V(u, v, w, x, y) + 2V(u^2, v, w, x, y) =
5[\bold 0] = \Delta^0_{\bold f} .$$

\vskip .2in

Hence, we find the Vogel sets all have correct dimension, and so the Vogel cycles are defined and  $\Delta^3_{\bold f} =
V(x,y)$,
$\Delta^2_{\bold f} = 2V(u^2 + v^2 + w^2, x, y)$,
$\Delta^1_{\bold f} = 2V(u^2 + v^2, w, x, y)$, and $\Delta^3_{\bold f} = 5[\bold 0]$.  

\vskip .4in

\noindent{\it Remark 2.18}. Suppose that all the Vogel cycles of $\bold f$ are defined and $k+1>d$. Consider the {\it
truncated
$d$-tuple}
$\bold f_{\operatorname{tr}}:=(f_{k+1-d},
\dots, f_k)$; we claim that, in a neighborhood of $V(\bold f)$,  $|V(\bold f)|= |V({\bold f_{\operatorname{tr}}})|$ and  both
$\bold f$ and
$\bold f_{\operatorname{tr}}$  will produce the same
$D^i$, $\Delta^i$, 
$\Pi^{i}$,
$\widetilde\Pi^{i}$, and $\widehat\Pi^{i}$ for all $i$ (all of them will be empty for $i<0$).

It is immediate from the definitions that $\widetilde\Pi^{i}_\bold f = \widetilde\Pi^{i}_{\bold f_{\operatorname{tr}}}$ and
$\widehat\Pi^{i}_\bold f = \widehat\Pi^{i}_{\bold f_{\operatorname{tr}}}$.  We would know that, near $V(\bold f)$,
$\Pi^{i}_\bold f =
\Pi^{i}_{\bold f_{\operatorname{tr}}}$ and, hence, that $D^i_\bold f = D^i_{\bold f_{\operatorname{tr}}}$ and
$\Delta^i_\bold f =
\Delta^i_{\bold f_{\operatorname{tr}}}$, if we could show that there is an equality of sets 
$|V(\bold f)| = |V({\bold f_{\operatorname{tr}}})|$. 

This is easy; by definition of $\Pi^i_\bold f$, $|V(\bold f_{\operatorname{tr}})| = |V(\bold f)|\cup |\Pi^0_\bold f|$. As we
are assuming that $\Pi^0_\bold f$ is $0$-dimensional (and, of course, has no components contained in $V(\bold f)$), there is
a neighborhood of $V(\bold f)$ in which $|V(\bold f)| = |V({\bold f_{\operatorname{tr}}})|$.

\vskip .2in

Suppose that all the Vogel cycles of $\bold f$ are defined and $k+1<d$. Consider the   {\it extended $d$-tuple} $\bold
f_{\operatorname{ex}}:=(f_0, \dots, f_0, \dots, f_k)$ (where there are $d-k$ occurrences of
$f_0$); clearly, $|V(\bold f)|= |V({\bold f_{\operatorname{ex}}})|$, and $\bold f$ and $\bold f_{\operatorname{ex}}$  will
produce the same $D^i$, $\Delta^i$, $\Pi^{i}$, $\widetilde\Pi^{i}$, and $\widehat\Pi^{i}$ for all $i$ (all of them will be
empty for $i<d-(k+1)$).

\vskip .3in

 Looking at the two cases above, we see that, if all the Vogel cycles are defined, the whole theory remains unchanged if we
assume that
$k+1 = d$, i.e., if we assume that our tuple $\bold f$ contains as exactly as many functions as the dimension of the
underlying space $X$.

\vskip .4in

We now prove a theorem which gives the basic relation between Vogel cycles and the blow-up. In fact, we show that the Vogel
cycles are representatives of the {\it Segre classes}, as defined in [{\bf Fu}, \S 4.2].  In the generic case, this is
Theorem 3.3 of [{\bf Gas1}], and is also proved in Lemma 2.2 of [{\bf G-G}].  However, we are interested in cases which may
not be quite so generic.

\vskip .4in

\noindent{\bf Theorem 2.19}. {\it Let $X$ be an irreducible analytic subset of
an analytic manifold $\Cal U$, let $\pi:\operatorname{Bl}_{\bold f}X\rightarrow X$
denote the blow-up of
$X$ along $\bold f$ (see Example 1.6), and let $E_{\bold f}$ denote the corresponding exceptional divisor.

\vskip .1in

If $E_{\bold f}$ properly intersects $\Cal U\times\Bbb P^m\times\{\bold 0\}$ in $\Cal U\times\Bbb P^k$ for all $m$, then 

\vskip .1in

\noindent i)\hskip .2in the Vogel cycles of $\bold f$ are defined;

\vskip .1in

\noindent ii)\hskip .17in there exists a neighborhood $\Omega$ of $V(\bold f)$ such that, for all $m$,
$\operatorname{Bl}_{\bold f}X$ intersects $\nomathbreak{\Omega\times\Bbb P^m\times\{\bold 0\}}$\linebreak\hbox{}\hskip .3in 
properly in
$\Omega\times\Bbb P^k$; and

\vskip .1in

\noindent iii)\hskip .14in  inside $\Omega$, for
all $i$,
$$\widehat\Pi^{i+1}_{\bold f} = \pi_*(\operatorname{Bl}_{\bold f}X\cdot (\Cal U\times \Bbb P^{i+k+1-d}\times\{\bold 0\}))$$ and
$$\Delta^i_{\bold f} = \pi_*(E_{\bold f}\cdot (\Cal U\times \Bbb P^{i+k+1-d}\times\{\bold 0\})),$$ where the intersection takes
place in
$\Cal U\times \Bbb P^k$ and $\pi_*$ denotes the proper push-forward.

\vskip .1in

Moreover, for all $\bold p\in X$, there exists an open neighborhood $\Cal U$ of $\bold p$ in $\Cal U$ such that, for a generic
linear reorganization,
$\tilde{\bold f}$, of
$\bold f$, 
$E_{\tilde{\bold f}}$ properly intersects
$\Cal U\times\Bbb P^m\times\{\bold 0\}$ inside
$\Cal U\times \Bbb P^k$ for all $m$. In the algebraic category, we may produce such generic linear reorganizations globally,
i.e., such that $E_{\tilde{\bold f}}$ properly intersects
$\Cal U\times\Bbb P^m\times\{\bold 0\}$ inside
$\Cal U\times \Bbb P^k$ for all $m$.

}

\vskip .3in

\noindent{\it Proof}. We show the last two statements first. As in 2.11, the reason that we can only make local statements
in the analytic case is because we must worry about analytic sets having an infinite number of irreducible components. For
all $\bold p\in X$, $\pi^{-1}(\bold p)$ is compact, and so, any analytic set can have only a finite number of irreducible
components which meet $\pi^{-1}(\bold p)$. In the algebraic setting, we know that we have a finite number of irreducible
components globally. For notational ease, we assume in the following paragraph, in the analytic case, that $\Cal U$ is rechosen as
small as necessary at each stage so that $\Cal U\times\Bbb P^k$ contains a finite number of analytic components (of any specified
analytic set) which intersect $\pi^{-1}(\bold p)$; this will mean that we will write $\Cal U$ in place of the open
neighborhood $\Cal U$, which appears in the statement of the theorem.

Now, as each point in each component of
$E_{\bold f}$ cannot have all of its homogeneous coordinates equal to zero, for each component $\nu$ of $E_{\bold f}$, there
exists a homogeneous coordinate
$w_{k(\nu)}$ such that $V(w_{k(\nu)})$ properly intersects $\nu$. Therefore, for generic $(a_{0,0}, \dots, a_{0,k})\in\Bbb
C^{k+1}$, the linear form $\widetilde w_k := a_{0,0}w_0 + \dots+ a_{0,k}w_k$ is such that $V(\widetilde w_k)$ contains no
component of
$E_{\bold f}$. We continue in this manner; for generic $(a_{1,0}, \dots, a_{1,k})\in\Bbb C^{k+1}$, the linear form
$\widetilde w_{k-1} := a_{1,0}w_0 + \dots+ a_{1,k}w_k$ is such that $V(\widetilde w_{k-1})$ contains no component of
$E_{\bold f}\cap V(\widetilde w_k)$. Continuing, we produce a generic linear reorganization,
$\widetilde{\bold w}$, of $\bold w$ such that, for all $m$, $E_{\bold f}$ properly intersects $V(\widetilde w_{m+1},
\dots,
\widetilde w_k)$ inside $\Cal U\times \Bbb P^k$. This proves the last two claims of the theorem.

\vskip .3in

We now prove i), ii), and iii) of the theorem.

\vskip .2in

We use $[w_0:\dots:w_k]$ as homogeneous coordinates on $\Bbb P^k$. Let
$\eta:\operatorname{Bl}_{\bold f}X\rightarrow
\Bbb P^k$ denote the restriction of the projection. Until the end of the proof, we shall simply write $f_j$ in place of
$f_j\circ\pi$; no confusion will arise, since it is clear that we must mean $f_j\circ\pi$ when the domain is contained in
$\operatorname{Bl}_{\bold f}X$.

Certainly, $\pi^{-1}$ induces an isomorphism from $\Pi^{i+1}_{\bold f}-V(\bold f)$ to
$$\eta^{-1}(\Bbb P^{i+k+1-d}\times\{\bold 0\})-E_\bold f \ =\ \operatorname{Bl}_\bold f X\cap (\Cal U\times\Bbb
P^{i+k+1-d}\times\{\bold 0\})-E_\bold f.$$  Hence, $\Pi^{i+1}_{\bold f}$ is purely $(i+1)$-dimensional if and only if 
$$
\overline{\operatorname{Bl}_\bold f X\cap (\Cal U\times\Bbb P^{i+k+1-d}\times\{\bold 0\})-E_\bold f} 
$$ is purely $(i+1)$-dimensional. But, every component of $\operatorname{Bl}_\bold f X\cap (\Cal U\times\Bbb
P^{i+k+1-d}\times\{\bold 0\})$ has dimension at least $i+1$, while -- by hypothesis -- $E_\bold f\cap (\Cal U\times\Bbb
P^{i+k+1-d}\times\{\bold 0\})$ is purely
$i$-dimensional. Thus, 
$$
\overline{\operatorname{Bl}_\bold f X\cap (\Cal U\times\Bbb P^{i+k+1-d}\times\{\bold 0\})-E_\bold f} =
\operatorname{Bl}_\bold f X\cap (\Cal U\times\Bbb P^{i+k+1-d}\times\{\bold 0\}),
$$ and every component has dimension at least $i+1$.  As $E_\bold f$ is locally defined in $\operatorname{Bl}_\bold f X$ by
a single equation and $E_\bold f\cap (\Cal U\times\Bbb P^{i+k+1-d}\times\{\bold 0\})$ is purely $i$-dimensional, it follows that
$\operatorname{Bl}_\bold f X\cap (\Cal U\times\Bbb P^{i+k+1-d}\times\{\bold 0\})$ is purely $(i+1)$-dimensional, for all
$i$, at all points which lie in $E_\bold f$.  This proves ii) from the statement of the theorem, and proves that
$\Pi^{i+1}_{\bold f}$ is purely $(i+1)$-dimensional, for all
$i$, at all points of $V(\bold f)$, and so the Vogel cycles are defined. This proves i).

\vskip .1in

Note that the Dimensionality Lemma and the above paragraphs imply that, in a neighborhood of any point $\bold p\in V(\bold
f)$,
$$
\operatorname{Bl}_\bold f X\cap V(w_{i+k+2-d}, \dots, w_k) = \overline{\operatorname{Bl}_\bold f X\cap V(w_{i+k+2-d},
\dots, w_k)-V(f_{i+k+1-d})}. \tag{*}
$$ 

\vskip .2in

Let $\bold p$ be a point in $V(\bold f)$. As the Vogel cycles are defined, there exists a neighborhood of $\bold p$ such that
$X-V(\bold f)$, $V(f_k)-V(\bold f)$, \dots, $V(f_{i+k+2-d})-V(\bold f)$ all intersect properly and $\pi$ induces an
isomorphism 
$$\big[\operatorname{Bl}_\bold f X-E\big]\cdot \big[V(f_k)-E\big]\cdot\dots\cdot \big[V(f_{i+k+2-d})-E\big]\cong$$
 $$\big[X-V(\bold f)\big]\cdot \big[V(f_k)-V(\bold f)\big]\cdot\dots\cdot \big[V(f_{i+k+2-d})-V(\bold f)\big].$$ By the
Dimensionality Lemma, no component of this intersection is contained in $V(f_{i+k+1-d})$, and so we conclude that
$\widehat\Pi^{i+1}_{\bold f}$ is equal to 
$$\pi_*\left(\overline{\big[\operatorname{Bl}_\bold f X-V(f_{i+k+1-d})\big]\cdot
\big[V(f_k)-V(f_{i+k+1-d})\big]\cdot\dots\cdot\big[V(f_{i+k+2-d})-V(f_{i+k+1-d})\big]}\right).$$ We claim that this implies
the first equality of the theorem:
$$
\widehat\Pi^{i+1}_{\bold f}=\pi_*\big(\operatorname{Bl}_\bold f X\cdot V(w_{i+k+2-d}, \dots, w_k)\big),\tag{$\dagger$}
$$ in a neighborhood of any point in $V(\bold f)$.

\vskip .3in

To see this, note that $\operatorname{Bl}_\bold f X-V(f_{i+k+1-d})\subseteq\{w_{i+k+1-d}\neq 0\}$. On the open set,
$\Cal W\subseteq \Cal U\times \Bbb P^k$, where
$f_{i+k+1-d}\neq 0$ and $w_{i+k+1-d}\neq 0$, there is an equality of schemes 
$$\operatorname{Bl}_\bold f X = V\left(\frac{f_j}{f_{i+k+1-d}}-\frac{w_j}{w_{i+k+1-d}}\right)_{j\neq i+k+1-d}.$$ At points of
$\Cal W$, $\dsize
\left\{\frac{f_j}{f_{i+k+1-d}}-\frac{w_j}{w_{i+k+1-d}}\right\}_{j\neq i+k+1-d}$ is easily seen to be a regular sequence.
Therefore, on $\Cal W$, the cycle $\left[\operatorname{Bl}_\bold f X\right]$ is equal to the intersection product of the
cycles 
$$\left[V\left(\frac{f_j}{f_{i+k+1-d}}-\frac{w_j}{w_{i+k+1-d}}\right)\right]_{j\neq i+k+1-d}.$$ Moreover, on $\Cal W$, for
$j\geqslant i+k+2-d$,
$$
\left[V\left(\frac{f_j}{f_{i+k+1-d}}-\frac{w_j}{w_{i+k+1-d}}\right)\right]\cdot \big[V(f_j)\big] \ =\
\left[V\left(\frac{f_j}{f_{i+k+1-d}}-\frac{w_j}{w_{i+k+1-d}}, \ f_j\right)\right] \ =$$

$$\big[V(f_j, w_j)\big]\ =\ \big[V(f_j)\big]\cdot\big[V(w_j)\big] \ =$$ 

$$
\left[V\left(\frac{f_j}{f_{i+k+1-d}}-\frac{w_j}{w_{i+k+1-d}}, \ w_j\right)\right] \ =\
\left[V\left(\frac{f_j}{f_{i+k+1-d}}-\frac{w_j}{w_{i+k+1-d}}\right)\right]\cdot \big[V(w_j)\big].$$ Hence, on $\Cal W$, 
$$
\left[\operatorname{Bl}_\bold f X\right]\cdot \big[V(f_k)\big]\cdot\dots\cdot\big[V(f_{i+k+2-d})\big]\ =
\left[\operatorname{Bl}_\bold f X\right]\cdot \big[V(w_k)\big]\cdot\dots\cdot\big[V(w_{i+k+2-d})\big]\ =
$$
$$
\left[\operatorname{Bl}_\bold f X\right]\cdot \big[V(w_{i+k+2-d}, \dots, w_k)\big],
$$ and so $(\dagger)$ follows from our previous paragraphs and $(*)$.

Now, by definition, $\Delta^i_\bold f+\widehat\Pi^i_\bold f = \widehat\Pi^{i+1}_\bold f\cdot V(f_{i+k+1-d})$. Applying
$(\dagger)$ and the push-forward formula (see Appendix A.14) -- which we may use since $V(f_{i+k+1-d}\circ\pi)$ properly
intersects 
$\operatorname{Bl}_\bold f X\cap V(w_{i+k+2-d}, \dots, w_k)$ by (*) -- we conclude that 
$$
\Delta^i_\bold f+\widehat\Pi^i_\bold f = \pi_*\big(V(f_{i+k+1-d}\circ\pi)\cdot\operatorname{Bl}_\bold f X\cdot V(w_{i+k+2-d},
\dots, w_k)\big).
$$ By the Dimensionality Lemma, $\Delta^i_\bold f$ consists of those components of the proper push-forward which are
contained in
$V(f)$.  Hence, we will have proved the second equality of the theorem if we can show that the components of
$V(f_{i+k+1-d}\circ\pi)\cdot\operatorname{Bl}_\bold f X\cdot V(w_{i+k+2-d},
\dots, w_k)$ which are contained in $E_\bold f$ are equal to $E_\bold f\cdot V(w_{i+k+2-d},
\dots, w_k)$.

On the open set where $w_{i+k+1-d}\neq 0$, $E_\bold f$ is defined to be
$V(f_{i+k+1-d}\circ\pi)\cdot\operatorname{Bl}_\bold f X$. Thus, it is enough to show that
$V(f_{i+k+1-d}\circ\pi)\cdot\operatorname{Bl}_\bold f X\cdot V(w_{i+k+2-d},
\dots, w_k)$ has no components contained in
$E_\bold f$ which are also contained in $V(w_{i+k+1-d})$. However, by hypothesis, $V(w_{i+k+1-d},
\dots, w_k)$ properly intersects $E_\bold f$, and so every component of $E_\bold f\cap V(w_{i+k+1-d}, w_{i+k+2-d},
\dots, w_k)$ has dimension $i-1$. As every component of
$V(f_{i+k+1-d}\circ\pi)\cdot\operatorname{Bl}_\bold f X\cdot V(w_{i+k+2-d},
\dots, w_k)$ has dimension at least $i$, we are finished.
\qed

\vskip .3in

\noindent{\it Remark 2.20}. Note that the proof of 2.19 shows that, for each $i$, if $E_{\bold f}$ properly
intersects $\Cal U\times \Bbb P^{i+k+1-d}\times\{\bold 0\}$ in $\Cal U\times\Bbb P^k$, then $\Pi^{i+1}_{\bold f}$ is purely
$(i+1)$-dimensional near $V(\bold f)$ -- the point being that we do {\bf not} need to assume that we have proper intersections for
{\bf all} $i$.

\vskip .4in

The following corollary follows immediately from Theorem 2.19.

\vskip .4in

\noindent{\bf Corollary 2.21 (The Segre-Vogel Relation)}. {\it Let $X$ be an analytic subset of
an analytic manifold $\Cal U$, and let $\pi:\Cal U\times\Bbb P^k\rightarrow \Cal U$
denote the projection. Assume that each $V_l$ appearing in $M$ is $d$-dimensional. For each $V_l$ appearing in $M$, consider 
${\operatorname{Bl}}_{\bold f}V_l\subseteq V_l\times\Bbb P^k\subseteq\Cal U\times\Bbb P^k$, and let
$E^l_{\bold f}$ denote the corresponding exceptional divisor. Let ${\operatorname{Bl}}_{\bold f}M:=\sum_l m_l[{\operatorname{Bl}}_{\bold
f}V_l]$ and $E_{\bold f}(M):=\sum_l m_l[E^l_{\bold f}]$.

\vskip .1in

If $|E_{\bold f}(M)|$ properly intersects $\Cal U\times\Bbb P^m\times\{\bold 0\}$ in $\Cal U\times\Bbb P^k$ for all $m$, then 

\vskip .1in

\noindent i)\hskip .2in the Vogel cycles of $\bold f$ with respect to $M$ are defined;

\vskip .1in

\noindent ii)\hskip .17in there exists a neighborhood $\Omega$ of $|M|\cap V(\bold f)$ such that, for all $m$,
$|\operatorname{Bl}_{\bold f}M|$ intersects $\Omega\times\Bbb P^m\times\{\bold 0\}$\linebreak\hbox{}\hskip .3in 
properly in
$\Omega\times\Bbb P^k$; and

\vskip .1in

\noindent iii)\hskip .14in  inside $\Omega$, for
all $i$,
$$\widehat\Pi^{i+1}_{\bold f}(M) = \pi_*(\operatorname{Bl}_{\bold f}M\cdot (\Cal U\times \Bbb P^{i+k+1-d}\times\{\bold 0\}))$$ and
$$\Delta^i_{\bold f}(M) = \pi_*(E_{\bold f}(M)\cdot (\Cal U\times \Bbb P^{i+k+1-d}\times\{\bold 0\})),$$ where the intersection
takes place in
$\Cal U\times \Bbb P^k$ and $\pi_*$ denotes the proper push-forward.

\vskip .1in

Moreover, for all $\bold p\in |M|\cap X$, there exists an open neighborhood $\Cal W$ of $\bold p$ in $\Cal U$ such that, for a
generic linear reorganization,
$\tilde{\bold f}$, of $\bold f$, 
$|E_{\tilde{\bold f}}(M)|$ properly intersects
$\Cal U\times\Bbb P^m\times\{\bold 0\}$ inside
$\Cal U\times \Bbb P^k$ for all $m$. In the algebraic category, we may produce such generic linear reorganizations globally,
i.e., such that $|E_{\tilde{\bold f}}(M)|$ properly intersects
$\Cal U\times\Bbb P^m\times\{\bold 0\}$ inside
$\Cal U\times \Bbb P^k$ for all $m$.
}

\vskip .3in

\noindent{\bf Definition 2.22}. We call a generic linear reorganization of $\bold f$, such as appears in Corollary 2.21, a
{\it Vogel reorganization of $\bold f$ with respect to $M$}. 

A generic linear reorganization of $\bold f$ which is both agreeable and Vogel is called {\it unifying}. 

\vskip .3in

\noindent{\it Remark 2.23}. Theorem 3.3 of [{\bf Gas1}] actually shows that, by replacing $\bold f$ by a generic linear
transformation applied to $\bold f$, one obtains a unifying $\tilde\bold f$; the point being that the linear transformation
is actually {\it generic}, not just generic in the IPZ topology. However, as one can see in the proof of 2.19, proving that
one can use an IPZ-generic transformation to obtain a suitable $\tilde\bold f$ is quite trivial, and is actually what one
uses in examples. 

\vskip .4in

\noindent   \S3.  {\bf The L\^e-Iomdine-Vogel Formulas}

\vskip .2in 

As in the previous section, $X$ will denote an analytic space of dimension $d$ contained in an analytic manifold $\Cal U$, $\bold f
:= (f_0, \dots, f_k)$  will be an ordered $(k+1)$-tuple of elements of $\Cal O_{{}_X}$, and $M=\sum_l m_l [V_l]$ will be an
analytic cycle in $X$ such that $\pm M>0$. 

We wish to examine the effect on the Vogel cycles of adding scalar multiples of a large power of a new
function
$g: X\rightarrow\Bbb C$ to
$f_0$. The formulas that we derive are a powerful tool for inductive proofs.

Throughout most of this section, we will be making the assumption that the Vogel cycles of $\bold f$ have correct dimension;
as discussed in Remark 2.18, this means that we may as well assume that the number of elements of  $\bold f$ is exactly $d$.
Therefore, we will find it convenient to let $n := d-1$, and then write that the dimension of $X$ is
$n+1$ and that
$\bold f = (f_0, \dots, f_n)$. Moreover, as all of our results will concern gap and Vogel {\bf cycles}, the contributions
from various irreducible components of $M$ will simply add, and so -- for simplicity -- we will make the assumption that $X$
is irreducible (though not necessarily reduced) and prove most results in the case where $M=[X]$.

Since we will be assuming that the gap sets have correct dimension, $\Delta^0_\bold f$ will be purely $0$-dimensional, and
for any $\bold p\in X$, we write $\big(\Delta^0_\bold f\big)_\bold p$ for the coefficient (possibly zero) of $\bold p$
appearing in the cycle $\Delta^0_\bold f$.

\vskip .3in

The following lemma relates the Vogel cycles of $(f_0, \dots, f_n)$ to the Vogel cycles of $(f_1,
\dots, f_n, g)$, where $g$ is a new function. We think of this as relating the Vogel cycles of $\bold f$ to the Vogel cycles
of
$\bold f$ restricted to $V(g)$ -- the elimination of $f_0$ corresponds to the drop in dimension of the ambient space. As we
shall see later, this ``restriction'' lemma is an essential step in proving the L\^e-Iomdine-Vogel (LIV) formulas.

\vskip .3in

\noindent{\bf Lemma 3.1 (The Restriction Lemma)}. {\it Let $X$ be an irreducible analytic space of dimension $n+1$, and let
$\bold f:= (f_0,
\dots, f_n)\in
\big(\Cal O_X\big)^{n+1}$. Let $g\in\Cal O_X$, let $\bold h:=(f_1, \dots, f_n, g)$,  and let $\bold p\in V(\bold f, g)$.  

\vskip .1in

\noindent i)\hskip .22in Suppose that $\Pi^1_\bold f$ is purely
$1$-dimensional at $\bold p$. Then, $\Pi^1_\bold f$ properly intersects $V(g)$ at $\bold p$ if and only if $V(\bold h) =
V(\bold f, g)$ as germs of sets at
$\bold p$.

\vskip .2in

\noindent ii)\hskip .2in Suppose that the Vogel sets of $\bold f$ have correct dimension at $\bold p$, that $V(\bold h) =
V(\bold f, g)$ as germs of sets at
$\bold p$, and that $V(g)$ properly intersects $D^i_\bold f$ at $\bold p$ for all $i\geqslant 1$. 

Then, ${\operatorname{dim}}_\bold p V(\bold h) = \big({\operatorname{dim}}_\bold p V(\bold f)\big) - 1$ provided that 
${\operatorname{dim}}_\bold p V(\bold f)\geqslant 1$, $V(g)$ properly intersects
$\widehat\Pi^i_\bold f$ at $\bold p$ for all
$i$, the Vogel sets of
$\bold h$  have correct dimension at $\bold p$, and, for all $i$ such that $1\leqslant i\leqslant n$, there are equalities
of germs of cycles at
$\bold p$ given by
$$
\widehat\Pi^i_\bold h = \widehat\Pi^{i+1}_\bold f\cdot V(g)\ \text{ and }\ \Delta^i_\bold h = \Delta^{i+1}_\bold f\cdot V(g).
$$ In addition, when $i=0$, we have the following equality of germs of cycles
$$
\Delta^0_\bold h = \big(\widehat\Pi^1_\bold f\cdot V(g)\big) + \big(\Delta^1_\bold f\cdot V(g)\big) .
$$

}

\vskip .3in

\noindent{\it Proof}. 

\vskip .1in

\noindent Proof of i):\hskip .2in As germs of sets at $\bold p$, 
$$V(f_1, \dots, f_n, g) = \big(\Pi^1_\bold f\cup V(\bold f)\big)\cap V(g) = \big(\Pi^1_\bold f\cap V(g)\big)\cup V(\bold f,
g).
$$ Since $\Pi^1_\bold f$ is purely
$1$-dimensional at $\bold p$, that $\Pi^1_\bold f$ properly intersects $V(g)$ at $\bold p$ is equivalent to
$\Pi^1_\bold f\cap V(g)$ being empty or equal to $\{\bold p\}$. As $\bold p\in V(\bold f, g)$, we have proved i).

\vskip .1in

\noindent Proof of ii):\hskip .2in By 2.4, $V(\bold f) = \bigcup D^i_\bold f$. As the Vogel cycles have correct dimension
and those of dimension at least one are properly intersected by $V(g)$ at
$\bold p$, we conclude that  ${\operatorname{dim}}_\bold p V(\bold h)$ equals $\big({\operatorname{dim}}_\bold p V(\bold
f)\big) - 1$ provided that 
${\operatorname{dim}}_\bold p V(\bold f)\geqslant 1$.

To see that $V(g)$ properly intersects
$\widehat\Pi^i_\bold f$ at $\bold p$, we work solely with germs of sets at $\bold p$. For $i\leqslant 0$,
$\Pi^i_\bold f =
\emptyset$, and so there is nothing to prove. We now proceed with a proof by contradiction. Let $m$ be the smallest $i$ such
that
$\Pi^i_\bold f$ does not properly intersect $V(g)$ at $\bold p$. Note that i) implies that $m\geqslant 2$, and we have that
${\operatorname{dim}}_\bold p \Pi^m_\bold f\cap V(g) = m$. Thus, ${\operatorname{dim}}_\bold p \Pi^m_\bold f\cap
V(f_{m-1})\cap V(g)
\geqslant m$. However, $\Pi^m_\bold f\cap V(f_{m-1}) = \Pi^{m-1}_\bold f\cup D^{m-1}_\bold f$, and so we would have to have
that either ${\operatorname{dim}}_\bold p \Pi^{m-1}_\bold f\cap V(g) \geqslant m-1$ or ${\operatorname{dim}}_\bold p
D^{m-1}_\bold f\cap V(g) \geqslant m-1$; the first possibility is excluded by definition of $m$, and the second possibility
is excluded by hypothesis. Thus, we have shown that $V(g)$ properly intersects $\widehat\Pi^i_\bold f$ at
$\bold p$ for all $i$.

To show that the Vogel sets of
$\bold h$  have correct dimension at $\bold p$, we once again work on the level of germs of sets. By definition,
$\Pi^i_\bold h = V(f_{i+1},
\dots, f_n, g)\
\lnot\ V(\bold h)$. One of our assumptions is that $V(\bold h)= V(\bold f, g)$; hence, $\Pi^i_\bold h = V(f_{i+1},
\dots, f_n, g)\
\lnot\ V(\bold f, g)$. We apply 1.3.iii to obtain $\Pi^i_\bold h =
\big(\Pi^{i+1}_\bold f\cap V(g)\big)\
\lnot\ V(\bold f, g)$. However, $V(\bold f, g)$ contains no components of $\Pi^{i+1}_\bold f\cap V(g)$, for
$\Pi^{i+1}_\bold f\cap V(g)$ is purely $i$-dimensional, while -- as sets --
$\Pi^{i+1}_\bold f\cap V(\bold f)\cap V(g) = \big(\bigcup_{m\leqslant i}D^m_\bold f\big)\cap V(g)$, which has dimension less
than
$i$. Therefore, as germs of sets at $\bold p$, $\Pi^i_\bold h = \Pi^{i+1}_\bold f\cap V(g)$ and is purely
$i$-dimensional, and so the Vogel sets of $\bold h$ have correct dimension at $\bold p$.

We wish to see that, for $1\leqslant i\leqslant n$, $\widehat\Pi^i_\bold h = \widehat\Pi^{i+1}_\bold f\cdot V(g)$ at
$\bold p$. As we saw above, this equality holds for the underlying sets and neither set has a component contained in
$V(\bold f)$. Therefore, it is enough to show that the cycles $\widehat\Pi^i_\bold h$ and $\widehat\Pi^{i+1}_\bold f\cdot
V(g)$ are equal on $X-V(\bold f)$. Applying Remark 2.13, we find that both of these cycles on $X-V(\bold f)$ are given by
$V(f_{i+1})\cdot\dots\cdot V(f_n)\cdot V(g)$.

Finally, for $0\leqslant i\leqslant n-1$, $\Delta^i_\bold h = \widehat\Pi^{i+1}_\bold h\cdot V(f_{i+1}) -
\widehat\Pi^i_\bold h$. Thus, for $1\leqslant i\leqslant n-1$, 
$$
\Delta^i_\bold h\ =\ \widehat\Pi^{i+2}_\bold f\cdot V(g)\cdot V(f_{i+1}) - \widehat\Pi^{i+1}_\bold f\cdot V(g)\ =\ 
\big(\widehat\Pi^{i+1}_\bold f+\Delta^{i+1}_\bold f\big)\cdot V(g) - \widehat\Pi^{i+1}_\bold f\cdot V(g) =
\Delta^{i+1}_\bold f\cdot V(g).
$$

When $i=0$, we have 
$$\Delta^0_\bold h = \widehat\Pi^{1}_\bold h\cdot V(f_{1}) = \widehat\Pi^{2}_\bold f\cdot V(g)\cdot V(f_{1}) =
\big(\widehat\Pi^{1}_\bold f+\Delta^{1}_\bold f\big)\cdot V(g).
$$

When $i=n$, $\Delta^n_\bold h\ =\ \widehat\Pi^{n+1}_\bold h\cdot V(g) - \widehat\Pi^n_\bold h\ =\
\widehat\Pi^{n+1}_\bold h\cdot V(g) - \widehat\Pi^{n+1}_\bold f\cdot V(g)$. We need to show that
$\widehat\Pi^{n+1}_\bold h - \widehat\Pi^{n+1}_\bold f =
\Delta^{n+1}_\bold f$. 

If $\bold f\equiv 0$, then $\widehat\Pi^{n+1}_\bold f = 0$, $\Delta^{n+1}_\bold f = [X]$, and -- as $V(g)$ properly
intersects $\Delta^{n+1}_\bold f$ -- we conclude that $\bold h\not\equiv 0$ and so $\widehat\Pi^{n+1}_\bold h = [X]$. Thus,
if
$\bold f\equiv 0$, $\widehat\Pi^{n+1}_\bold h - \widehat\Pi^{n+1}_\bold f =
\Delta^{n+1}_\bold f$. If $\bold f\not\equiv 0$, then $\widehat\Pi^{n+1}_\bold f = [X]$, $\Delta^{n+1}_\bold f = 0$, and --
as $V(g)$ properly intersects $\widehat\Pi^{n+1}_\bold f$ -- we conclude that $\bold h\not\equiv 0$ and so
$\widehat\Pi^{n+1}_\bold h = [X]$. Thus, if
$\bold f\not\equiv 0$, $\widehat\Pi^{n+1}_\bold h - \widehat\Pi^{n+1}_\bold f =
\Delta^{n+1}_\bold f$.\qed

\vskip .4in

\noindent{\bf Definition 3.2}. Suppose that $\bold p\in\widehat\Pi^1_\bold f\cap V(g)$ and that
${\operatorname{dim}}_\bold p
\widehat\Pi^1_\bold f = 1$. Let $\eta$ be an irreducible component (with its reduced structure) of $\widehat\Pi^1_\bold f$
which passes through $\bold p$.

If $\eta\cap V(g)$ is zero-dimensional at $\bold p$, then we define {\it the gap ratio of $\eta$ at $\bold p$} (for
$\bold f$ with respect to $g$) to be the ratio of intersection numbers $\dsize\frac{\big(\eta\cdot V(f_{k+1-d})\big)_\bold
p}{\big(\eta\cdot V(g)\big)_\bold p}$.

\vskip .1in

If $\eta\cap V(g)$ is not zero-dimensional at $\bold p$ (i.e., if $\eta\subseteq V(g)$), then we define {\it the gap ratio of
$\eta$ at
$\bold p$} (for
$\bold f$ with respect to $g$) to be $0$.

\vskip .1in

A {\it gap ratio} (at $\bold p$ for $\bold f$ with respect to $g$) is any one of the gap ratios of any component of
$\widehat\Pi^1_\bold f\cap V(g)$ through $\bold p$.

\vskip .1in

If $\bold p\in V(g)$, but $\bold p\not\in\widehat\Pi^1_\bold f$, then we say that {\it all the gap ratios are zero}.

\vskip .1in

Finally, a {\it gap ratio at $\bold p$ for $\bold f$ with respect to $g$ and the cycle $M$} is a gap ratio (at $\bold p$ for $\bold
f$ with respect to $g$) of ${\bold f}_{|_{V_l}}$ for some $V_l$ appearing in $M$.

\vskip .4in

\noindent{\bf Lemma 3.3}. {\it Let $X$ be an irreducible analytic space of dimension $n+1$, let $\bold f:= (f_0, \dots,
f_n)\in
\big(\Cal O_X\big)^{n+1}$, let $g\in\Cal O_X$, and let $\bold p\in V(g)$. Let $a$ be a non-zero complex number, and let
$j\geqslant 1$ be an integer. 

\vskip .1in

If $j$ is greater than or equal to the maximum gap ratio at $\bold p$ for $\bold f$ with respect to $g$, then, for all but
(possibly) a finite number of complex $a$,

\vskip .2in

\noindent i)\hskip .24in  if $\widehat\Pi^1_\bold f$ is purely one-dimensional at $\bold p$, then $\widehat\Pi^1_\bold f$
properly intersects
$V(f_0 + ag^j)$ at
$\bold p$, and $\big(\Delta^0_\bold f\big)_\bold p =
\big(\widehat\Pi^1_\bold f\cdot  V(f_0 + ag^j)\big)_\bold p$.

\vskip .2in

Moreover, if we have the strict inequality that $j$ is greater than the maximum gap ratio at $\bold p$ for $\bold f$ with
respect to $g$, then i) holds for all non-zero $a$; in particular, this is the case if $j\geqslant 1+
\big(\Delta^0_\bold f\big)_\bold p$.

\vskip .2in

\noindent ii)\hskip .2in Suppose that $\Pi^1_\bold f$ is purely $1$-dimensional at $\bold p$, and that $\bold p\in V(\bold
f, g)$. Then, $\Pi^1_\bold f$ properly intersects
$V(f_0 + ag^j)$ at
$\bold p$ if and only if  there is an equality of germs of sets at
$\bold p$ given by
$V(f_1,
\dots, f_n, f_0+ag^j) = V(\bold f, g)$.

\vskip .2in

\noindent iii)\hskip .2in Suppose that $\bold p\in V(\bold f, g)$ and that, at $\bold p$,  there is an equality of germs of
sets given by
$$V(f_1,
\dots, f_n, f_0+ag^j) = V(\bold f, g),$$ the Vogel sets of $\bold f$ all have correct dimension, and that, for all 
$i\geqslant 1$, $V(g)$ properly intersects each $D^i_\bold f$.

 If
$1\leqslant i\leqslant n$, then, at
$\bold p$,
$\widehat\Pi^{i+1}_\bold f$ properly intersects
$V(f_0+ag^j)$, the Vogel sets of the $(n+1)$-tuple $(f_1, \dots, f_n, f_0+ag^j)$ have correct dimension, and there is an
equality of germs of cycles given by
$$
\widehat\Pi^i_{(f_1, \dots, f_n, f_0+ag^j)}\ =\ \widehat\Pi^{i+1}_\bold f\cdot V(f_0+ag^j).
$$ }

\vskip .3in

\noindent{\it Proof}.

\vskip .1in

\noindent i)\hskip .2in Assume that $\widehat\Pi^1_\bold f$ is purely one-dimensional at $\bold p$. Then, we may write the
cycle
$\big[\widehat\Pi^1_\bold f\big]$ as $\sum m_\nu[\nu]$, where each $\nu$ is a reduced, irreducible curve at $\bold p$. Let
$\alpha_\nu(t)$ denote a local parameterization of $\nu$ such that $\alpha_\nu(0)=\bold p$. Then, to show that
$\widehat\Pi^1_\bold f$ properly intersects $V(f_0 + ag^j)$ at $\bold p$, we need to show that, for all $\nu$,
$(f_0 + ag^j)_{|_{\alpha_\nu(t)}}\not\equiv 0$. To show that $\big(\Delta^0_\bold f\big)_\bold p =
\big(\widehat\Pi^1_\bold f\cdot  V(f_0 + ag^j)\big)_\bold p$, we need to show that $\big(\widehat\Pi^1_\bold f\cdot 
V(f_0)\big)_\bold p = \sum m_\nu \big([\nu]\cdot  V(f_0 + ag^j)\big)_\bold p$; calculating intersection numbers as in A.9 of
Appendix A, we find that what we need to show is that, for all $\nu$, ${\operatorname{mult}}_t f_0(\alpha_\nu (t)) =
{\operatorname{mult}}_t ((f_0 + ag^j)\circ\alpha_\nu)(t)$. Thus, we may prove both the proper intersection statement and the
intersection formula at the same time by proving this multiplicity statement.

	Clearly, ${\operatorname{mult}}_t ((f_0 + ag^j)\circ\alpha_\nu)(t)= \operatorname{min}\big\{{\operatorname{mult}}_t
(f_0\circ\alpha_\nu)(t),
\ {\operatorname{mult}}_t (g^j\circ\alpha_\nu)(t) \big\}$,  unless the lowest degree terms of $f_0(\alpha_\nu(t))$ and
$-a(g^j\circ\alpha_\nu)(t)$ are precisely equal. As
${\operatorname{mult}}_t (f_0\circ\alpha_\nu)(t) = \big([\nu]\cdot V(f_0)\big)_\bold p$ and ${\operatorname{mult}}_t
(g^j\circ\alpha_\nu)(t)= j\big([\nu]\cdot V(g)\big)_\bold p$, we conclude that ${\operatorname{mult}}_t ((f_0 +
ag^j)\circ\alpha_\nu)(t) = {\operatorname{mult}}_t (f_0\circ\alpha_\nu)(t)$ if $j$ is greater than the maximum gap ratio,
and that this equality holds when $j$ equals the maximum gap ratio except for the finite number of values of $a$ which would
cause cancellation of the lowest degree terms. This proves i).

\vskip .2in

\noindent ii)\hskip .2in This follows immediately by applying Lemma 3.1.i with the $g$ of the lemma replaced by
$f_0+ag^j$.
\vskip .2in

\noindent iii)\hskip .2in This follows immediately by applying Lemma 3.1.ii with the $g$ of the lemma replaced by
$f_0+ag^j$.\qed

\vskip .4in

\noindent{\bf Theorem 3.4 (The L\^e-Iomdine-Vogel formulas)}. {\it Suppose that each $V_l$ appearing in $M$ has dimension
$n+1$. Let $\bold f:= (f_0, \dots, f_n)\in
\big(\Cal O_X\big)^{n+1}$, let $g\in\Cal O_X$, and let $\bold p\in |M|\cap V(\bold f, g)$. Let $a$ be a non-zero complex number,
let
$j\geqslant 1$ be an integer, and let $\bold h := (f_1, \dots, f_n, f_0+ag^j)$.

Suppose that the Vogel cycles of $\bold f$ with respect to $M$ are defined at $\bold p$, and that $V(g)$ properly intersects each
of the Vogel cycles, $\Delta^i_\bold f(M)$, at
$\bold p$  for all $i\geqslant 1$.

\vskip .1in

If $j$ is greater than or equal to the maximum gap ratio at $\bold p$ for $\bold f$ with respect to $g$ and $M$, then for all but
(possibly) a finite number of complex $a$, in a neighborhood of $\bold p$:

\vskip .1in

\noindent i)\hskip .2in there is an equality of sets given by
$|M|\cap V(\bold h) = |M|\cap V(\bold f, g)$,

\vskip .1in

\noindent ii)\hskip .16in ${\operatorname{dim}_\bold p}(|M|\cap V(\bold h)) = \big({\operatorname{dim}_\bold p}\big(|M|\cap V(\bold
f)\big)\big) -1$ provided that
${\operatorname{dim}_\bold p}(|M|\cap V(\bold f))\geqslant 1$,

\vskip .1in

\noindent iii)\hskip .12in  the Vogel cycles of $\bold h$ with respect to $M$ exist at $\bold p$, and

\vskip .2in

\noindent iv)\hskip .14in $\Delta^0_\bold h(M) = \Delta^0_\bold f(M)+j\big(\Delta^1_\bold f(M)\cdot V(g)\big)$\hskip .1in and, for
$1\leqslant i\leqslant n-1$,\hskip .1in $\Delta^i_\bold h(M) = j\big(\Delta^{i+1}_\bold f(M)\cdot V(g)\big)$.

\vskip .3in

Moreover, if we have the strict inequality that $j$ is greater than the maximum gap ratio at $\bold p$ for $\bold f$ with
respect to $g$ and $M$, then these equalities hold for all non-zero $a$; in particular, this is the case if $j\geqslant 1+
{\operatorname{max}}_l\{\big(\Delta^0_{{\bold f}_{|_{V_l}}}\big)_\bold p\}$.

}

\vskip .3in

\noindent{\it Proof}.  The assumption that all $m_l$ have the same sign prevents cancellation of contributions from various
$V_l$; thus, the assumption that $V(g)$ properly intersects each  $\Delta^i_\bold f(M)$ implies that $V(g)$ properly intersects
each  $\Delta^i_{{\bold f}_{|_{V_l}}}$ for all $l$. Therefore, we are reduced to considering the case of Lemma 3.3, where $X$ is
irreducible and $M$ equals [X].

\vskip .1in

Now, the equality of sets in i) is precisely 3.3.ii; the statement concerning ${\operatorname{dim}_\bold
p}V(\bold h)$ follows from this equality of sets, combined with the facts that $V(\bold f) = \bigcup D^i_\bold f$ (see 2.4)
and that $V(g)$ properly intersects the non-zero-dimensional Vogel cycles of $V(\bold f)$.

\vskip .1in

Now, suppose that $0\leqslant i\leqslant n-1$. By definition, $\widehat\Pi^i_\bold h + \Delta^i_\bold h =
\widehat\Pi^{i+1}_\bold h\cdot V(f_{i+1})$. By 3.3.iii, this equals
$\widehat\Pi^{i+2}_\bold f\cdot V(f_0+ag^j)\cdot V(f_{i+1})$. By definition of the Vogel cycles, this equals
$\big(\widehat\Pi^{i+1}_\bold f+\Delta^{i+1}_\bold f\big)\cdot V(f_0+ag^j)$.  As $\big|\Delta^{i+1}_\bold f\big|\subseteq
V(f_0)$, $\Delta^{i+1}_\bold f\cdot V(f_0+ag^j) = \Delta^{i+1}_\bold f\cdot V(ag^j) = j\big(\Delta^{i+1}_\bold f\cdot
V(g)\big)$. Therefore, we have shown that 
$$
\widehat\Pi^i_\bold h + \Delta^i_\bold h = \big(\widehat\Pi^{i+1}_\bold f\cdot V(f_0+ag^j)\big) + j\big(\Delta^{i+1}_\bold
f\cdot V(g)\big).\tag{$\dagger$}
$$

\vskip .1in

If $i=0$, then $\widehat\Pi^i_\bold h = 0$, and the first equality of the theorem follows from $(\dagger)$ and 3.3.i.

\vskip .1in

If $1\leqslant i\leqslant n-1$, then 3.3.iii tells us that $\big(\widehat\Pi^{i+1}_\bold f\cdot V(f_0+ag^j)\big) =
\widehat\Pi^i_\bold h$; cancelling $\widehat\Pi^i_\bold h$ from each side of $(\dagger)$ yields the second equality of the
theorem.\qed

\vskip .4in

\noindent{\it Remark 3.5}. A principal use of the LIV formulas is in families;  one requires something about the
constancy of the Vogel cycles of $\bold f$ in the family, and the LIV formulas imply the constancy of the Vogel cycles of a tuple
of function with a smaller zero locus.

However, it is possible to use these formulas ``in reverse'' -- to calculate the Vogel cycles of $\bold h_{(a, j)}:=(f_1,
\dots, f_n, f_0+ag^j)$ and have them tell us about the Vogel cycles of $(f_0, \dots, f_n)$. The difficulty of applying the
LIV formulas in this manner is that it is not so easy to know when $j$ is greater than or equal to the maximum gap ratio. 
We discuss this problem below, using the notation from the theorem.

\vskip .1in

Suppose that the Vogel cycles of $\bold f$ are defined at $\bold p$, and that $V(g)$ properly intersects each of the Vogel
cycles, $\Delta^i_\bold f$, at
$\bold p$  for all $i\geqslant 1$.  Assume that, in a neighborhood of $\bold p$, there is an equality of sets given by
$V(\bold h_{(a, j)}) = V(\bold f, g)$ (we are still assuming that $a\neq 0$). 

By assuming that $V(g)$ properly intersects $\Delta^i_\bold f$ for $i\geqslant 1$, we are assuming that we can calculate the
Vogel sets of $\bold f$ in dimensions one and higher. While it would be nice to be able to proceed without this assumption,
there seems to be no way to avoid it. Notice that, if we could calculate
$(\Delta^0_\bold f)_\bold p$, then we would know that the LIV formulas hold for
$j>(\Delta^0_\bold f)_\bold p$. However, $(\Delta^0_\bold f)_\bold p$ is typically more difficult to calculate than
$(\Delta^1_\bold f\cdot V(g))_\bold p$. So, we will assume that we can also calculate the intersection number
$(\Delta^1_\bold f\cdot V(g))_\bold p$, and then consider the problem of how can one tell when $j$ is large enough for the
LIV formulas to hold using data gathered from
$\Delta^0_{{\bold h}_{(a, j)}}$ and $(\Delta^1_\bold f\cdot V(g))_\bold p$.

\vskip .2in

\noindent Our best answer is that:

\vskip .2in

 if $j\ > \ \big(\Delta^0_{{\bold h}_{(a, j)}}\big)_\bold p - j\big(\Delta^1_\bold f\big)_\bold p$,  then the LIV formulas
hold, and so $\big(\Delta^0_\bold f\big)_\bold p = \big(\Delta^0_{{\bold h}_{(a, j)}}\big)_\bold p - j\big(\Delta^1_\bold
f\big)_\bold p$. 

\vskip .3in

To see this, note that the proof of 3.4 shows that $\big(\Delta^0_{{\bold h}_{(a, j)}}\big)_\bold p - j\big(\Delta^1_\bold
f\big)_\bold p =
\big(\widehat\Pi^1_\bold f\cdot V(f_0+ag^j)\big)_\bold p$. We claim that, if $j >\big(\widehat\Pi^1_\bold f\cdot
V(f_0+ag^j)\big)_\bold p$, then $j>(\Delta^0_\bold f)_\bold p = \big(\widehat\Pi^1_\bold f\cdot V(f_0)\big)_\bold p$ and so
the LIV formulas hold. This is easy; calculating intersection numbers as in the proof of 3.3, 
$$
\big(\widehat\Pi^1_\bold f\cdot V(f_0+ag^j)\big)_\bold p \ \geqslant \ \operatorname{min}\big\{\big(\widehat\Pi^1_\bold
f\cdot V(f_0)\big)_\bold p , \  j\big(\widehat\Pi^1_\bold f\cdot V(g)\big)_\bold p \big\}.
$$ The desired conclusion follows. 

\vskip .3in

One might hope that if $\big(\Delta^0_{{\bold h}_{(a, j+1)}}\big)_\bold p - \big(\Delta^0_{{\bold h}_{(a, j)}}\big)_\bold p =
\big(\Delta^1_\bold f\big)_\bold p$  (which would be true if the LIV formulas held),  then one could, in fact, conclude that
the LIV formulas {\bf do} hold. Unfortunately, the situation is slightly more complicated than this.

  Let us call
$(a, j)$ an {\it exceptional pair} if there exists a component $\nu$ of $\widehat\Pi^1_\bold f$ at $\bold p$ such that 
$$
\big(\nu\cdot V(f_0+ag^j)\big)_\bold p \ \neq \ \operatorname{min}\big\{\big(\nu\cdot V(f_0)\big)_\bold p,\ j\big(\nu\cdot
V(g)\big)_\bold p\big\}.
$$ Looking at the proofs of 3.3 and 3.4, it is easy to see that, if $(a, j)$ is {\bf not} an exceptional pair, then 
$\big(\Delta^0_{{\bold h}_{(a, j+1)}}\big)_\bold p - \big(\Delta^0_{{\bold h}_{(a, j)}}\big)_\bold p \geqslant
\big(\Delta^1_\bold f\big)_\bold p$ with equality if and only if $j$ is greater than or equal to the maximum polar ratio.
Hence, if it were not for the existence of exceptional pairs, one could simply make a table of values of
$\big(\Delta^0_{{\bold h}_{(a, j)}}\big)_\bold p$ for fixed $a$ and increasing $j$, and when a difference between successive
entries is exactly $\big(\Delta^1_\bold f\big)_\bold p$, one would have identified the maximum polar ratio and would know
that the LIV formulas hold beyond that value for $j$.

On the other hand, if $(a, j)$ {\bf is} an exceptional pair, it is quite possible that $\big(\Delta^0_{{\bold h}_{(a,
j+1)}}\big)_\bold p -
\big(\Delta^0_{{\bold h}_{(a, j)}}\big)_\bold p =
\big(\Delta^1_\bold f\big)_\bold p$ and still $j$ is smaller than the maximum polar ratio. Of course, this can only happen
once for each possible exceptional pair, and the number of exceptional pairs is certainly no more than the number of
components of $\widehat\Pi^1_\bold f$ through $\bold p$. Thus, if we know the number  of components of
$\widehat\Pi^1_\bold f$ through $\bold p$, call this number $c$, and we make a table of values of
$\big(\Delta^0_{{\bold h}_{(a, j)}}\big)_\bold p$, once we see a difference between successive values equalling 
$\big(\Delta^1_\bold f\big)_\bold p$ more than $c$ times, we know that $j$ is high enough for the LIV formulas to hold.

Alternatively, and only pseudo-rigorously, if one selects the constant $a$ ``randomly'', then $a$ will not be part of an
exceptional pair and so, $\big(\Delta^0_{{\bold h}_{(a, j+1)}}\big)_\bold p - \big(\Delta^0_{{\bold h}_{(a, j)}}\big)_\bold
p \geqslant
\big(\Delta^1_\bold f\big)_\bold p$ with equality if and only if $j$ is greater than or equal to the maximum polar ratio.
This approach is particularly well-suited for computer calculation.

\vskip .4in

The following lemma is related to 3.3 and 3.4 and will be of use to us later.

\vskip .2in

\noindent{\bf Lemma 3.6 }. {\it Suppose that each $V_l$ appearing in $M$ has dimension
$n+1$, let $\bold f:= (f_0,
\dots, f_n)\in
\big(\Cal O_X\big)^{n+1}$, let $\bold p\in |M|\cap V(\bold f)$, and suppose that the Vogel cycles of $\bold f$ with respect to $M$ at $\bold
p$ exist. 

Let $a$ be a non-zero complex number, and let
$j\geqslant 1$ be an integer. Let $\pi$ denote the projection from
$\Bbb C\times X$ to
$X$, and let
$w$ denote the projection from
$\Bbb C\times X$ to $\Bbb C$. 

Then, the Vogel sets of $\bold h:=(w^j,  f_1\circ\pi, \dots, f_n\circ\pi, f_0\circ\pi+aw^j)$ with respect to $\Bbb C\times M$ have correct
dimension at
$(0, \bold p)$, for all $i\leqslant n+1$, $\Bbb C\times\widehat\Pi^i_\bold f(M)$ properly intersects
$V(f_0\circ\pi+aw^j)$, and there is an equality of germs of cycles at $(0, \bold p)$ given by
$$
\widehat\Pi^i_{\bold h}(\Bbb C\times M)\ = \ \big(\Bbb C\times\widehat\Pi^i_\bold f(M)\big)\ \cdot\ V(f_0\circ\pi+aw^j).
$$ 

}

\vskip .3in

\noindent{\it Proof}. As usual, we instantly reduce ourselves to the case where $X$ is irreducible and $M=[X]$.

\vskip .1in

That $\Bbb
C\times\widehat\Pi^i_\bold f$ properly intersects
$V(f_0\circ\pi+aw^j)$ is obvious.

\vskip .1in

First, note that $V(\bold h) = \{0\}\times V(\bold f)$. Suppose that $1\leqslant i\leqslant n+1$. Then,
$$
\Pi^i_\bold h = V(f_i\circ\pi, \dots, f_n\circ\pi, f_0\circ\pi+aw^j)\ \lnot\ \big(\{0\}\times V(\bold f)\big).
$$ We have $V(f_i\circ\pi, \dots, f_n\circ\pi) = \big(\Bbb C\times\Pi^i_\bold f\big)\cup V(\bold f\circ\pi)$, and
$V(\bold f\circ\pi)\cap V(f_0\circ\pi+aw^j)\subseteq \{0\}\times V(\bold f)$.  Applying 1.3.iii, we find that 
$$
\Pi^i_\bold h = \big(\Bbb C\times\Pi^i_\bold f\big)\cap V(f_0\circ\pi+aw^j)\ \lnot\ \big(\{0\}\times V(\bold f)\big).
$$ Now, near $\bold p$, $\big(\Bbb C\times\Pi^i_\bold f\big)\cap V(f_0\circ\pi+aw^j)$ is purely $i$-dimensional, and -- not
only does it have no components contained in $\{0\}\times V(\bold f)$ -- in fact, it has no components contained in
$\Bbb C\times V(\bold f)$; for, by 2.4, the set $\Pi^i_\bold f\cap V(\bold f)$ equals the union of all of the Vogel sets of
dimension less than or equal to $i-1$. Therefore, the set $\Pi^i_\bold h$ equals the set $\big(\Bbb C\times\Pi^i_\bold
f\big)\cap V(f_0\circ\pi+aw^j)$ and, hence, the Vogel sets of $\bold h$ have correct dimension at
$(0, \bold p)$. 

As we saw above, $\big|\widehat\Pi^i_\bold h\big| = \big|\Pi^i_\bold h\big|$ has no components contained in $\Bbb C\times
V(\bold f)$. Thus, to prove that the cycles $\widehat\Pi^i_{\bold h}$ and $\big(\Bbb C\times\widehat\Pi^i_\bold f\big)\
\cdot\ V(f_0\circ\pi+aw^j)$ are equal, it is enough to prove the equality on $\big(\Bbb C\times X\big)-\big(\Bbb C\times
V(\bold f)\big)$. Once again, we apply Remark 2.13 and find that both cycles are equal to 
$$ V(f_i\circ\pi)\cdot \dots \cdot V(f_n\circ\pi)\cdot V(f_0\circ\pi+aw^j)
$$ on $\big(\Bbb C\times X\big)-\big(\Bbb C\times V(\bold f)\big)$.
\qed

\vskip .4in

\noindent\S4. {\bf Concluding Remarks}

\vskip .2in

As we stated in the introduction, this paper describes the algebraic framework necessary to generalize our work on affine hypersurface
singularities to the case of functions on arbitrary analytic spaces. If, in fact, we were only interested in dealing with the case where the
underlying space was a local complete intersection (l.c.i.), then we could have avoided introducing the cycle $M$ into most of our results;
we wish to give a short explanation as to why this is the case.

\vskip .1in

Suppose that we have an analytic function $\tilde g:\Cal U\rightarrow\Bbb C$. In [{\bf M2}], we describe the
relationship between the vanishing cycles along $\tilde g$ (with respect to the constant sheaf on $\Cal U$)  and the exceptional divisor in
the blow-up of
$\Cal U$ along the Jacobian ideal; this relationship allowed us to derive many results concerning how the L\^e cycles and numbers controlled
the topological data provided by the vanishing cycles.

Now, let $g:=\tilde g_{|_{X}}$. We wish to produce a similar relationship between the exceptional divisor in the blow-up of {\bf some} ideal
and the vanishing cycles along
$g$ with respect to {\bf some} constructible complex of sheaves. What turns out to be crucial about the constant sheaf -- call it  $\bold
P^\bullet$ -- on
$\Cal U$ is that its characteristic cycle (see [{\bf K-S}]), $\operatorname{Ch}(\bold P^\bullet)$, is such that $\pm \operatorname{Ch}(\bold
P^\bullet)>0$; we know that $\bold P^\bullet$ has this property because the (shifted) constant sheaf on a manifold is a {\it perverse} sheaf
(see [{\bf BBD}]). Moreover, we show in [{\bf M1}] that the correct exceptional divisor to look at is in the blow-up of
$\operatorname{Ch}(\bold P^\bullet)$ along the ideal defining the image of differential
$d\tilde g$.

Unfortunately, if $X$ is not an l.c.i., then the constant sheaf on $X$ need not be perverse, and so there is no
guarantee that
$\pm\operatorname{Ch}(\bold P^\bullet)>0$, {\bf provided} that $\bold P^\bullet$ is the constant sheaf on $X$. However, one can take $\bold
P^\bullet$ to be the {\it perverse cohomology} ([{\bf BBD}]) of the constant sheaf on $X$, and then the results of this paper can be used to
show that the entire theory of L\^e cycles can be generalized to arbitrary underlying spaces.

\enddocument